\documentclass[11pt,a4paper]{article}

\usepackage{amsmath, amsfonts, amssymb, amsthm, amstext, amsbsy, bm, mathrsfs, esint, stmaryrd}

\usepackage[hmargin=1in, vmargin=0.75in]{geometry}
\usepackage{fancyhdr}
\usepackage{mathtools}

\usepackage{graphicx, color}
\usepackage{subcaption}

\usepackage{url, hyperref}

\usepackage{algorithm, algpseudocode}

\usepackage{enumerate}
\usepackage{enumitem}
\usepackage{verbatim, pifont, prettyref}

\font\msbm=msbm10

\theoremstyle{plain}
\newtheorem{theorem}{Theorem}[section]
\newtheorem{lemma}[theorem]{Lemma}

\newtheorem{proposition}[theorem]{Proposition}

\newtheorem{remark}[theorem]{Remark}

\numberwithin{equation}{section}

\newcommand{\RR}{\mathbb{R}}

\newcommand{\be}{\boldsymbol{e}}

\newcommand{\bg}{\boldsymbol{g}}

\newcommand{\bv}{\boldsymbol{v}}

\newcommand{\bx}{\boldsymbol{x}}

\newcommand{\biota}{\boldsymbol{\iota}}

\newcommand{\balpha}{\boldsymbol{\alpha}}
\newcommand{\bbeta}{\boldsymbol{\beta}}

\newcommand{\BA}{\boldsymbol{A}}
\newcommand{\BB}{\boldsymbol{B}}

\newcommand{\BH}{\boldsymbol{H}}
\newcommand{\BI}{\boldsymbol{I}}

\newcommand{\BM}{\boldsymbol{M}}
\newcommand{\BK}{\boldsymbol{K}}
\newcommand{\BP}{\boldsymbol{P}}

\newcommand{\BS}{\boldsymbol{S}}

\newcommand{\BV}{\boldsymbol{V}}

\newcommand{\BX}{\boldsymbol{X}}
\newcommand{\BY}{\boldsymbol{Y}}

\newcommand{\BG}{\boldsymbol{G}}

\newcommand{\BPhi}{\boldsymbol{\Phi}}

\newcommand{\bzero}{\boldsymbol{0}}
\newcommand{\PL}{\mathrm{P}\text{\L}}
\DeclareMathOperator{\Tr}{Tr}

\DeclareMathOperator{\diag}{diag}

\DeclareMathOperator{\sign}{sign}

\newcommand{\I}{\boldsymbol{I}}

\newcommand{\lag}{\langle}
\newcommand{\rag}{\rangle}

\definecolor{xl}{RGB}{200,50,120}

\def\mathbb#1{\hbox{\msbm{#1}}}

\begin{document}

\title{\bf Multi-Step Performance Estimation for Gradient Descent under the Polyak–Łojasiewicz Condition }
\author{Wanli Hong\thanks{Center for Data Science, New York University Shanghai, Shanghai, China. S.L. and W.H. are (partially) financially supported by the National Key R\&D Program of China, Project Number 2021YFA1002800, National Natural Science Foundation of China (NSFC) Grant No.12571105 and No.12001372, Science and Technology Commission of Shanghai Grant No.24QA2706200 and No.24ZR1455300, Shanghai Municipal Education Commission (SMEC) via Grant 0920000112, and NYU Shanghai Boost Fund. W.H. is also supported by NYU Shanghai Ph.D. fellowship and acknowledges the NSF/NRT support.} \thanks{Center for Data Science, New York University.}, Shuyang Ling$^*$}
\date{\today}
\maketitle
\begin{abstract}
Gradient descent is a fundamental algorithm in machine learning and statistical inference. Performance estimation provides a systematic framework for obtaining sharp convergence bounds through semidefinite programming.
We study the worst-case performance of constant-stepsize gradient descent over the class of $L$-smooth functions satisfying the $\mu$-Polyak--\L{}ojasiewicz (P\L{}) inequality. We formulate the $N$-step objective-gap analysis as a semidefinite programming (SDP) relaxation based on pairwise smoothness constraints and pointwise P\L{} inequalities. Since an exact interpolation condition for smooth P\L{} functions is not currently available, the resulting SDP provides a valid upper bound on the true worst-case performance. We solve a reduced SDP retaining only consecutive-pair smoothness constraints analytically for every horizon $N$ and every stepsize $t\in(0,2/L)$. The solution exhibits three stepsize regimes. For $0<t\leq 1/L$ and $1/L<t\leq \sqrt{3}/L$, the reduced and full pairwise SDPs have the same optimal value, and the one-step SDP contraction factor tensorizes to a geometric $N$-step rate. For $\sqrt{3}/L<t<2/L$, the reduced-SDP value is nongeometric in the horizon: an initial transient phase is followed by eventual geometric decay. The resulting asymptotic contraction factor is tight for the full pairwise SDP and strictly improves the bound obtained by recursively applying the one-step estimate. Explicit primal and dual feasible solutions certify the results, with numerical SDP experiments guiding their construction. Our results show that analyzing several iterations together can give a better convergence bound than repeatedly applying the sharp one-step SDP bound in the Regime III.
\end{abstract}

\section{Introduction}\label{s:intro}

Gradient descent (GD) is a fundamental algorithm in machine learning and statistical inference. Understanding its worst-case convergence remains a central question in optimization.
The performance estimation problem (PEP), introduced by Drori and Teboulle~\cite{DT14}, formulates worst-case analysis as an optimization problem over the objective function and the oracle data generated by an algorithm. For a fixed horizon, the relevant data are the sampled triples $(\bx_i,\bg_i,f_i)$, where $\bg_i=\nabla f(\bx_i)$ and $f_i=f(\bx_i)$. An exact interpolation theorem characterizes when a given set of finite data arise from a function in the target class. By using the interpolation condition, the worst-case analysis can be formulated into a semidefinite program (SDP) that involves the Gram-matrix of the gradients.
This principle applies to the PEP analysis for smooth convex and strongly convex functions~\cite{THG17,THG17b} and enables both numerical exploration and analytical worst-case certificates.

We consider the problem of minimizing a differentiable function
$f:\RR^d\to\RR$ with finite optimal value
\[
f^{\star}:=\inf_{\bx\in\RR^d}f(\bx)>-\infty,
\]
using gradient descent with a constant stepsize $t>0$,
\[
\bx_{k+1}=\bx_k-t\nabla f(\bx_k), \qquad k=0,1,\cdots
\]
Our focus is the class of $L$-smooth functions satisfying the
Polyak--\L{}ojasiewicz (P\L{}) inequality
\begin{equation}\label{eq:PL}
    \frac12\|\nabla f(\bx)\|^2 \geq\mu\bigl(f(\bx)-f^{\star}\bigr), \qquad \bx\in\RR^d.
\end{equation}

The gradient inequality was first found in Polyak~\cite{P63} and
\L{}ojasiewicz~\cite{Loj63}, and it has become a standard
assumption in optimization and machine learning~\cite{KNS16,LZB22}. It guarantees global optimality of stationary points without requiring convexity or isolated minimizers.

P\L-type conditions arise in  low-rank matrix recovery and factorization~\cite{LLSW19,SL16} and neural-network training~\cite{KNS16,DZPS19,LZB22}. In these settings, the P\L-type conditions and smoothness properties can yield linear convergence despite nonconvexity. This motivates understanding the precise convergence guarantees available under the smooth P\L~assumption alone.

For $0<\mu<L$, denote the smooth P\L~class by
\[
{\cal F}_{\mu,L}^{\PL}
:=\left\{ f:\RR^d\to\RR\ \middle|\
\begin{array}{l}
f^{\star}>-\infty,\quad \nabla f\text{ is }L\text{-Lipschitz},~
f\text{ satisfies }\eqref{eq:PL}
\end{array}
\right\}.
\]
For $0<t<2/L$, the descent lemma and~\eqref{eq:PL} give
\[
\begin{aligned}
    f(\bx_{k+1})-f^{\star}
    &\leq f(\bx_k)-f^{\star}
        -t\left(1-\frac{Lt}{2}\right)\|\nabla f(\bx_k)\|^2\\
    &\leq\bigl(1-\mu t(2-Lt)\bigr)
        \bigl(f(\bx_k)-f^{\star}\bigr).
\end{aligned}
\]
Therefore, the standard linear convergence bound is
\[
\begin{aligned}
 f(\bx_N)-f^{\star} \leq p_{\rm van}(t)^N \bigl(f(\bx_0)-f^{\star}\bigr),~~~p_{\rm van}(t):=1-\mu t(2-Lt).
\end{aligned}
\]
To investigate whether this decay factor, or a sharper one-step factor, captures the worst-case behavior over several iterations, we define
\begin{equation}\label{def:pep_v0}
\mathcal R_N(t)
:=\sup\left\{ \frac{f(\bx_N)-f^{\star}}{f(\bx_0)-f^{\star}}~ \Big| 
    \begin{array}{l}
    f\in{\cal F}_{\mu,L}^{\PL},\quad f(\bx_0)>f^{\star},\\
    \bx_{k+1}=\bx_k-t\nabla f(\bx_k),\quad 0\leq k\leq N-1
    \end{array}
\right\}.
\end{equation}

\subsection{From one-step to multi-step performance estimation}

For smooth P\L~functions, exact interpolation remains the main obstacle. Every finite sample satisfies the interpolation inequalities for ${\cal F}_{-L,L}$ (defined in Section~\ref{s:prelim}) and the pointwise P\L~constraints $\|\bg_i\|^2\geq2\mu(f_i-f^{\star})$, but these conditions are not sufficient for interpolation by a globally smooth P\L~function. We use them to define the \emph{full pairwise SDP}~\eqref{def:sdp_pep}, with normalized value $\overline{\cal R}_N(t)$, and the further \emph{consecutive-pair SDP}~\eqref{primal:reSDP}, with normalized value $\widehat{\cal R}_N(t)$. The latter retains only the smoothness inequalities between consecutive iterates. Section~\ref{s:prelim} derives these relaxations, which satisfy
\begin{equation}\label{eq:three-pep-values}
    \mathcal R_N(t)
    \leq\overline{\cal R}_N(t)
    \leq\widehat{\cal R}_N(t).
\end{equation}
Constructing the primal--dual certificates and verifying the optimality condition establish exact SDP values. To infer an exact functional worst-case value ${\cal R}_N(t)$, one additionally needs to establish the interpolation condition for ${\cal F}_{\mu,L}^{\PL}$ or construct a matching smooth P\L{} function that attains the SDP bound, which remains unknown in general.

The closest prior P\L-PEP result is due to Abbaszadehpeivasti,
de Klerk, and Zamani~\cite{AKZ23}. Their earlier work~\cite{AKZ22} analyzes constant step GD for nonconvex $L$-smooth functions using semidefinite performance estimation and establishes exactness in certain stepsize regimes. Under the additional P\L~condition,~\cite{AKZ23} solves the one-step SDP relaxation and obtains the factor
\[
  p_{\rm ADZ}(t):=\widehat{\cal R}_1(t) =\overline{\cal R}_1(t).
\]
Its expression changes across three stepsize ranges, separated by $1/L$ and $\sqrt{3}/L$. Applying the one-step inequality recursively yields
\begin{equation}\label{eq:iterated-one-step}
    \mathcal R_N(t)\leq p_{\rm ADZ}(t)^N.
\end{equation}
However, sharpness for one iteration does not imply sharpness for a fixed multi-step horizon: extremal one-step data need not concatenate into a jointly feasible trajectory. This is the distinction addressed by our multi-step analysis.

At the level of the SDP relaxations studied here, we aim to answer this question:
\begin{center}
\emph{Does the sharp one-step SDP factor tensorize over multiple iterations?}
\end{center}
Here tensorization means that the $N$th power of the one-step factor is the exact value of the corresponding $N$-step SDP.

\subsection{Our contributions}

We determine the multi-step behavior of constant-stepsize GD within the semidefinite framework for the performance estimation problem. The results identify when the sharp one-step SDP factor tensorizes. Throughout, the {\em exactness} refers to the specified SDP relaxations~\eqref{def:sdp_pep} and~\eqref{primal:reSDP}, not automatically to the functional PEP~\eqref{def:pep_v0}.

\begin{enumerate}
\item {\bf Exact solution of the consecutive-pair SDP.}
For every $N\geq1$ and $0<t<2/L$, we determine
$\widehat{\cal R}_N(t)$ analytically, i.e., the SDP relaxation~\eqref{primal:reSDP} of the PEP based on all the consecutive steps. In Regimes~I and~II, $0<t\leq1/L$ and $1/L<t\leq\sqrt{3}/L$, respectively, its value $\widehat{{\cal R}}_N(t)$ is $p(t)^N$ with $p(t) = p_{\rm ADZ}(t)$. In Regime~III, $\sqrt{3}/L<t<2/L$, a
nongeometric transient is followed by a geometric phase:
\begin{equation}\label{eq:intro-regime-three-form}
\widehat{\cal R}_N(t)
=\begin{cases}
    c^{(N)}(t), & N\leq N_0(t),\\
    c^{(N_0)}(t)\,p(t)^{N-N_0(t)}, & N>N_0(t).
\end{cases}
\end{equation}
Our main result Theorem~\ref{thm:main} gives the explicit expressions for $c^{(N)}(t)$, $p(t)$, and the transition horizon $N_0(t)$.

\item {\bf Finite-horizon agreement with the full pairwise SDP.}
In Regimes~I and~II, our consecutive-pair optimizer satisfies every nonconsecutive smoothness constraint. Therefore, it holds that
\[
\overline{\cal R}_N(t) =\widehat{\cal R}_N(t) =p(t)^N,  \qquad N\geq1,
\]
so the one-step factor of~\cite{AKZ23} tensorizes exactly at the full-SDP level. The same feasibility property holds during the transient phase of Regime~III, giving
\[
\overline{\cal R}_N(t) =\widehat{\cal R}_N(t) =c^{(N)}(t)\leq c^{(1)}(t)^N= p_{\rm ADZ}(t)^N,  \qquad 1\leq N\leq N_0(t).
\]
Therefore, the multi-step relaxation captures nongeometric dependence on the horizon that is lost when the one-step bound is applied recursively; full-pairwise feasibility establishes that no further improvement arises from the omitted constraints in this phase.

\item {\bf Asymptotic tightness in Regime~III.}
For $N>N_0(t)$, the particular mixed optimizer in our construction violates some nonconsecutive constraints. We therefore do not establish finite-horizon equality between the consecutive-pair and full pairwise SDPs in this phase. Nevertheless, for every fixed $\sqrt{3}/L<t<2/L$, we prove
\[
\lim_{N\to\infty}\overline{\cal R}_N(t)^{1/N} =\lim_{N\to\infty}\widehat{\cal R}_N(t)^{1/N} =p(t)< p_{\rm ADZ}(t).
\]
Thus the reduced relaxation recovers the exact asymptotic geometric factor of the full pairwise SDP, strictly improving the asymptotic bound obtained by recursively applying the sharp one-step estimate.
\end{enumerate}

Explicit primal and dual feasible solutions certify the results.
Numerical SDP experiments guide the construction of the certificates, which are then verified analytically for arbitrary horizons. The full-SDP finite-horizon value in the geometric phase of Regime~III is not determined by this analysis, and all SDP values remain certified upper bounds on $\mathcal R_N(t)$ in the absence of a matching smooth P\L~interpolator.

\subsection{Related work}

Beyond convergence guarantees, the P\L~condition is closely related to error bounds~\cite{LT93} and quadratic growth~\cite{BI95,M00}. For $\mathcal C^2$ functions near a possibly non-isolated minimizer, Rebjock and Boumal~\cite{RB25} establish equivalences between the local P\L~condition and several such regularity properties. Under stronger global smoothness assumptions, Boumal, Criscitiello, and Rebjock~\cite{BCR26} show that smooth, globally P\L~functions admit a nonlinear least-squares representation. These results clarify the geometry of the P\L~class, complementing its role in convergence analysis.

Related interpolation frameworks cover weakly convex functions with bounded subgradients~\cite{AJ26}, smooth and strongly convex sets~\cite{AB24}, and linear operators~\cite{NJF24}. For smooth convex GD, sharp convergence analyses include the PEP bounds of Drori and Teboulle~\cite{DT14} and the elementary tight bounds of Teboulle and Vaisbourd~\cite{TV23}. PEPit supports numerical exploration of such formulations~\cite{BCF22}, while Das Gupta, Van Parys, and Ryu~\cite{DVR24} develop a branch-and-bound PEP framework for constructing methods with certifiably optimal worst-case guarantees. These developments illustrate the broader use of PEP for analysis and algorithm design; our focus is the analytical solution for constant step gradient descent under the P\L~condition.

Rubbens, Hendrickx, and Taylor~\cite{RHT26} develop a constructive one-point strengthening procedure for finite descriptions of function classes. Applied to smooth functions satisfying a \L{}ojasiewicz condition, it yields stronger necessary inequalities. For the smooth P\L~class, these conditions admit an SDP representation and can yield
strictly tighter numerical one-step bounds than the relaxation
in~\cite{AKZ23}, but do not constitute a complete exact interpolation theorem. Their work provides a strengthened numerical formulation, rather than a closed-form solution of the associated multi-step GD PEP. Our contribution is complementary: we solve a simpler relaxation analytically for every horizon and constant stepsize in $(0,2/L)$. Thus, our aim is to characterize its dependence on the horizon and
stepsize, not to obtain the smallest numerical upper bound available from strengthened constraints.

Kong, Jiang, and He~\cite{KJH26} also use PEP to derive explicit one-step rates for proximal-gradient methods under the P\L~inequality and a proximal-residual variant. Their analysis concerns composite optimization, while ours specializes to smooth GD and examines the joint constraints over multiple iterations. This comparison separates the extension to a broader algorithmic setting from the multi-step
question considered here.

In terms of methodology, Rotaru, Glineur, and Patrinos~\cite{RGP26} derive exact multi-step worst-case rates for constant-stepsize GD over classes with prescribed lower and upper curvature bounds. Their criterion is the minimum gradient norm along the trajectory, whereas we
study the terminal objective gap under the additional P\L~condition. Exact interpolation allows them to construct matching worst-case functions; our formulation remains a relaxation of the smooth P\L~class. Their analysis also highlights the role of nonconsecutive constraints for large stepsizes, with pairs at distance at most two sufficing in their setting. We similarly investigate when consecutive-pair constraints determine the full pairwise SDP value and
when our reduced-SDP optimizer fails to satisfy the omitted constraints.

Beyond constant stepsizes, nonconstant schedules can improve GD's worst-case guarantees over smooth convex and strongly convex classes~\cite{AP25a,AP25b,G24,GSW25}. These analyses are intrinsically multi-step: individual long steps need not decrease the objective, so the behavior of an entire schedule matters. Extending the present smooth P\L~analysis to prescribed or optimized schedules is a natural
direction, but lies outside our constant stepsize setting.
Other P\L~results address different notions of sharpness. Under additional local $\mathcal C^2$ regularity, Kassing and
Kruse~\cite{KK26} obtain asymptotic rates near a P\L~minimizer matching those of strongly convex quadratics. Their result concerns local asymptotic behavior, rather than global finite-horizon bounds under $L$-smoothness. At the oracle-complexity level,~\cite{YFL23} establishes an
$\Omega((L/\mu)\log(1/\varepsilon))$ lower bound on first-order oracle calls for general smooth P\L~optimization, showing that GD is optimal up to numerical constants. Our focus is instead on explicit finite-horizon SDP constants and their dependence on the fixed stepsize.

The remainder of the paper is organized as follows. Section~\ref{s:prelim} introduces the SDP relaxations and their dual formulations. Section~\ref{s:main} presents the main results and constructs the primal--dual certificates for the
three stepsize regimes. The appendices collect the supporting lemmas and their proofs.

\subsection{Notation}

Bold uppercase and lowercase letters denote matrices and vectors respectively. The entries of $\BX$ and $\bx$ are written as $X_{i,j}$ and $x_i$; a bold subscripted symbol such as $\BX_{i,j}$ instead denotes a matrix indexed by $(i,j)$. The notation $\BX_{i:j,~i':j'}$ selects rows $i,\ldots,j$ and columns $i',\ldots,j'$. The superscripts in parentheses indicate parameter dependence, as in $c^{(N)}$ and $\BX^{(N)}$, and are omitted when unambiguous.

For vectors, $\|\cdot\|$ denotes the Euclidean norm; for matrices, it denotes the operator norm. We write $\|\BX\|_*$ and $\|\BX\|_F$ for the nuclear and Frobenius norms, and $\lag\BX,\BY\rag:=\Tr(\BX^\top\BY)$ for the Frobenius inner product of matrices of the same size. For symmetric $\BX$, the relations $\BX\succeq0$ and $\BX\succ 0$ denote positive semidefiniteness and positive definiteness. Finally, $\{\be_k^{(N)}\}_{k=0}^N$ is the standard basis in $\RR^{N+1}$, $\I_N$ the $N\times N$ identity matrix, $\bzero_{n,m}$ the
$n\times m$ zero matrix, and $\sign(z):=z/|z|$ for $z\neq0$.

\subsection{Disclosure of AI assistance.}
We used ChatGPT (OpenAI) for language editing,
manuscript organization, algebraic checks, and proof simplification. For several supporting lemmas, we had already developed proofs based on recursion and induction on the time horizon, establishing a base case and extending the argument to general $N$. AI assistance contributed substantial simplifications of these arguments, notably by suggesting and drafting the weighted rank-one decompositions of the dual slack matrices in Lemmas~\ref{R1:L3} and~\ref{r3:LemmaS}, and the Kac--Murdock--Szeg\H{o} kernel argument for Gram-matrix positive definiteness in Lemma~\ref{R2:L3}. These structural arguments replaced lengthy recursive calculations with more direct and transparent proofs. Additional assistance included streamlining the interpolation and scalar induction arguments and identifying algebraic, indexing, and boundary-case errors. We retain full responsibility for the verification and correctness of all mathematical claims, proofs, computations, and references, and for the final content of the manuscript.

\section{Preliminaries and PEP}
\label{s:prelim}

This section derives the semidefinite relaxations and states the dual of the consecutive-pair relaxation. Throughout, we assume $0<\mu<L$, $0<t<2/L$, and $N\geq 1$. For a gradient-descent trajectory $\{\bx_k\}_{k=0}^N$, we write
\[
\bg_k:=\nabla f(\bx_k), \qquad f_k:=f(\bx_k), \qquad  \delta_k:=f_k-f^{\star}.
\]

\subsection{Smooth interpolation and the functional PEP}

To relate ${\cal F}_{\mu,L}^{\PL}$ to standard curvature classes, let
\[
{\cal F}_{m,L} :=\left\{
    f:\RR^d\to\RR\ \middle|\
    \begin{array}{l}
    \nabla f\text{ is }L\text{-Lipschitz},\\
    \bx\mapsto f(\bx)-\dfrac m2\|\bx\|^2\text{ is convex}
    \end{array}
\right\},
\qquad -L\leq m\leq L.
\]
Thus ${\cal F}_{0,L}$ and ${\cal F}_{\mu,L}$ are the smooth convex
and smooth strongly convex classes, respectively, whereas
${\cal F}_{-L,L}$ is the class of $L$-smooth and $L$-weakly convex functions. We have
\[
    {\cal F}_{\mu,L} \subseteq{\cal F}_{\mu,L}^{\PL} \subseteq{\cal F}_{-L,L}.
\]
The first inclusion follows from strong convexity and the second from
$L$-smoothness. In particular, it holds that
\[
 f(\bx)+\frac L2\|\bx\|^2 \in{\cal F}_{0,2L}, \qquad f\in{\cal F}_{\mu,L}^{\PL}.
\]
These inclusions allow smooth interpolation inequalities to be used in
the performance estimation formulation below, while distinguishing the
P\L~class from the strongly convex class.

We first recall the exact interpolation theorem for the class
$\mathcal F_{m,L}$.

\begin{theorem}[\bf Smooth interpolation for ${\cal F}_{m,L}$]
\label{thm:FmL-interpolation}
A finite set of triples $\{(\bx_i,\bg_i,f_i)\}_{i\in\mathcal I}$ is $\mathcal F_{m,L}$-interpolable with $m < L$ if and only if
\begin{equation}
\label{eq:FmL-interpolation}
\begin{aligned}
 f_i-f_j \geq  \lag \bg_j,\bx_i-\bx_j\rag +\frac{1}{2L}\|\bg_i-\bg_j\|^2 +\frac{m}{2(1-m/L)} \left\|\bx_i-\bx_j-\frac{1}{L}(\bg_i-\bg_j)\right\|^2
\end{aligned}
\end{equation}
for every $i,j\in\mathcal I$; see~\cite[Theorem~4]{THG17} and~\cite[Lemma 3.9]{THG17b}.
\end{theorem}

Since $\mathcal{F}_{\mu,L}^{\PL}\subseteq\mathcal F_{-L,L}$,
Theorem~\ref{thm:FmL-interpolation} with $m=-L$ implies that every finite sample from a function in $\mathcal{F}_{\mu,L}^{\PL}$ satisfies
\begin{equation}
\label{eq:PLsmooth}
\frac{1}{2L}\|\bg_i-\bg_j\|^2 -\frac{L}{4} \left\|\bx_i-\bx_j-\frac{1}{L}(\bg_i-\bg_j)\right\|^2 \leq f_i-f_j-\lag\bg_j,\bx_i-\bx_j\rag
\end{equation}
for every ordered pair $(i,j)$. Moreover, the P\L{} condition gives
\begin{equation}\label{eq:pointwise-PL}
    \delta_k\leq \frac{1}{2\mu}\|\bg_k\|^2.
\end{equation}
The conditions~\eqref{eq:PLsmooth}--\eqref{eq:pointwise-PL} are necessary but not sufficient. Thus, the finite-dimensional programs below are outer relaxations of the functional PEP.

For a fixed initial gap $\Delta>0$, the functional $N$-step PEP can be written as
\begin{equation}
\label{def:functional-pep}
\tag{$\mathrm{PEP}_N$}
\begin{aligned}
\sup_{f, \bx_0,\ldots,\bx_N, \ell}
\quad & \ell \\
{\rm s.t.}\quad
& \bx_{k+1}=\bx_k-t\nabla f(\bx_k),~~ 0\leq k\leq N-1,\\
& f(\bx_k)-f^{\star}\geq \ell,~~0\leq k\leq N,\\
& f(\bx_0)-f^{\star}\leq \Delta,\\
& f\in\mathcal{F}_{\mu,L}^{\PL}.
\end{aligned}
\end{equation}
With the notation introduced in Section~\ref{s:intro}, it holds that
\begin{equation}
\label{eq:functional-pep-value}
\operatorname{val}(\mathrm{PEP}_N) = \Delta {\cal R}_N(t).
\end{equation}
The normalized value is independent of $\Delta$ by homogeneity.

\subsection{Full and consecutive-pair SDP relaxations}

Every $L$-smooth function with finite infimum also satisfies
\[
\delta_k\geq \frac{1}{2L}\|\bg_k\|^2.
\]
Indeed, the descent lemma applied at
$\bx_k-L^{-1}\bg_k$ yields
$f^{\star}\leq f_k-(2L)^{-1}\|\bg_k\|^2$.

We now eliminate the iterates from~\eqref{eq:PLsmooth} and define
all pairwise coefficient matrices directly. Let
$\be_k^{(N)}\in\mathbb R^{N+1}$ be the $k$-th standard basis vector,
indexed by $k=0,\ldots,N$, and let
$\BG^{(N)}\in\mathbb S^{N+1}$ be the Gram matrix of
$\bg_0,\ldots,\bg_N$, with
\begin{equation}\label{eq:Gram-matrix}
 G_{i,j}^{(N)}:=\lag\bg_i,\bg_j\rag.
\end{equation}
For $0\leq i\neq j\leq N$, set $a:=\min\{i,j\}$ and $b:=\max\{i,j\}$. To simplify the notation, we write $\be_k=\be_k^{(N)}$ in the following display and define
\begin{equation}\label{def:Aij}
\begin{aligned}
 \BA_{i,j}^{(N)} & = \frac{1}{4L} (\be_i-\be_j)(\be_i-\be_j)^{\top}  \\
 &\quad + \frac{t \sign(j-i)}{4} \sum_{k=a}^{b-1} \Bigl[(\be_i+\be_j)\be_k^{\top} +\be_k(\be_i+\be_j)^{\top}\Bigr] -\frac{Lt^2}{4} \sum_{k=a}^{b-1}\sum_{\ell=a}^{b-1}\be_k\be_{\ell}^{\top}.
\end{aligned}
\end{equation}
Here $\sign(j-i)=1$ if $i<j$ and $-1$ if $i>j$. The matrix $\BA_{i,j}^{(N)}$ is symmetric and is supported on rows and columns $a,\ldots,b$; reversing the ordered pair changes only the sign of the middle term in~\eqref{def:Aij}. In particular, for every symmetric $\BG^{(N)}$, it holds that
\begin{equation}\label{eq:Aij-direct-inner-product}
\begin{aligned}
 \lag\BA_{i,j}^{(N)},\BG^{(N)}\rag
 &= \frac{G_{i,i}^{(N)}+G_{j,j}^{(N)}-2G_{i,j}^{(N)}}{4L}\\
 &\quad + \frac{t \sign(j-i)}{2} \sum_{k=a}^{b-1}\bigl(G_{i,k}^{(N)}+G_{j,k}^{(N)}\bigr) -\frac{Lt^2}{4} \sum_{k=a}^{b-1}\sum_{\ell=a}^{b-1}G_{k,\ell}^{(N)}.
\end{aligned}
\end{equation}
Substituting the gradient-descent  into
\eqref{eq:PLsmooth} and expanding the squared norms therefore
shows that its ordered-pair inequality is equivalent to
\[
 \lag\BA_{i,j}^{(N)},\BG^{(N)}\rag \leq\delta_i^{(N)}-\delta_j^{(N)},  \qquad 0\leq i\neq j\leq N.
 \]

The full pairwise SDP is
\begin{equation}\label{def:sdp_pep}\tag{$\overline{\mathrm P}_N$}
\begin{aligned}
\Delta \overline{{\cal R}}_N(t) 
:= \max_{\ell^{(N)}, \{\delta_k^{(N)}\}, \BG^{(N)}} \quad & \ell^{(N)} \\
{\rm s.t.}\quad & \lag\BA_{i,j}^{(N)},\BG^{(N)}\rag \leq \delta_i^{(N)}-\delta_j^{(N)}, ~~ 0\leq i\neq j\leq N,\\
& \frac{1}{2L}G_{k,k}^{(N)} \leq \delta_k^{(N)} \leq \frac{1}{2\mu}G_{k,k}^{(N)}, ~~ 0\leq k\leq N,\\
& \delta_0^{(N)}\leq \Delta,\\
& \delta_k^{(N)}\geq \ell^{(N)},~~ 0\leq k\leq N,\\
& \BG^{(N)}\succeq 0.
\end{aligned}
\end{equation}
Here $\overline{\cal R}_N(t)$ is the normalized optimal value of the full pairwise SDP and no rank constraint is imposed on $\BG^{(N)}$.

We next retain only the consecutive ordered pairs $(k,k+1)$ and $(k+1,k)$. The consecutive-pair relaxation is
\begin{align}
\label{primal:reSDP}
\tag{$\widehat{\mathrm P}_N$}
\Delta \widehat{{\cal R}}_N(t)
:=
\max_{\ell^{(N)}, \{\delta_k^{(N)}\}, \BG^{(N)}}
\quad & \ell^{(N)} \\
{\rm s.t.}\quad
& \lag\BA_{k,k+1}^{(N)},\BG^{(N)}\rag
 \leq \delta_k^{(N)}-\delta_{k+1}^{(N)}, ~~0\leq k\leq N-1,
 \qquad\tag{P1}\label{PP1} \\
& \lag\BA_{k+1,k}^{(N)},\BG^{(N)}\rag
 \leq \delta_{k+1}^{(N)}-\delta_k^{(N)},~~0\leq k\leq N-1,
 \qquad\tag{P2}\label{PP2}\\
& \frac{1}{2L}G_{k,k}^{(N)}
 \leq \delta_k^{(N)}
 \leq \frac{1}{2\mu}G_{k,k}^{(N)}, ~~ 0\leq k\leq N,
 \qquad\tag{P3}\label{PP3}\\
& \delta_0^{(N)}\leq \Delta,\qquad\tag{P4}\label{PP4}\\
& \delta_k^{(N)}\geq \ell^{(N)}, ~~0\leq k\leq N,
 \qquad\tag{P5}\label{PP5}\\
& \BG^{(N)}\succeq 0.  \qquad\tag{P6}\label{PP6}
\end{align}

The matrices $\BA_{k,k+1}^{(N)}$ and $\BA_{k+1,k}^{(N)}$ yield a simple form, i.e., they are basically equal to a $2\times 2$ matrix plus zero padding. To see this, we let
\[
 {\bf E}_k^{(N)} :=\begin{bmatrix}
 \be_k^{(N)} & \be_{k+1}^{(N)}
 \end{bmatrix}.
\]
A direct calculation gives
\[
\BA_{k,k+1}^{(N)} ={\bf E}_k^{(N)}\BA({\bf E}_k^{(N)})^{\top}, \qquad \BA_{k+1,k}^{(N)}  = {\bf E}_k^{(N)}\BA'({\bf E}_k^{(N)})^{\top},
\]
where
\begin{equation}
\label{def:AA}
\BA = \frac{1}{4L}
\begin{bmatrix}
2-s^2&s\\
s&1
\end{bmatrix},
\qquad
\BA' =\frac{1}{4L}
\begin{bmatrix}
-s^2-4s-2&-(s+2)\\
-(s+2)&1
\end{bmatrix},~~~s = Lt-1.
\end{equation}
Thus each consecutive coefficient matrix is supported on the corresponding $2\times2$ principal block and no separate support mapping is needed.

Since the first relaxation~\eqref{def:sdp_pep} uses only necessary smooth-P\L{} interpolation conditions and the second~\eqref{primal:reSDP} drops all nonconsecutive constraints, it holds that
\[
 {\cal R}_N(t) \leq \overline{{\cal R}}_N(t) \leq \widehat{{\cal R}}_N(t).
 \]
If an optimizer of $\widehat{\mathrm P}_N$ also satisfies every nonconsecutive constraint in $\overline{\mathrm P}_N$, then we can also show that  $\overline{{\cal R}}_N(t)=\widehat{{\cal R}}_N(t)$ because the primal and dual optimizers in~\eqref{primal:reSDP} still satisfy the KKT condition in~\eqref{def:sdp_pep}.

\paragraph{Dependence on the horizon.}
The superscript $(N)$ is retained whenever quantities from different horizons are compared. Let
\[
\biota_N := 
\begin{bmatrix}
\I_{N+1}\\
\bzero_{1,N+1}
\end{bmatrix}
 \in\RR^{(N+2)\times(N+1)}.
\]
For every $0\leq i\neq j\leq N$, we have
\begin{equation}\label{eq:A-horizon-embedding}
 \BA_{i,j}^{(N+1)}
 =\biota_N\BA_{i,j}^{(N)}\biota_N^{\top}.
\end{equation}
Thus the coefficient matrices are horizon-independent up to zero padding. By contrast, the optimal variables $\BG^{(N)}$, $\delta_k^{(N)}$, and $\lambda_k^{(N)}$ genuinely depend on $N$. 

\subsection{Dual of the consecutive-pair SDP}

We associate $\lambda_k^{(N)}\geq 0$ and
$\lambda_k^{\prime(N)}\geq 0$ with \textup{(P1)} and \textup{(P2)},
respectively. For \textup{(P3)}, $\alpha_k^{(N)}\geq 0$ is associated with
\[
 \delta_k^{(N)}-\frac{1}{2L}G_{k,k}^{(N)}\geq 0,
\]
whereas $\beta_k^{(N)}\geq 0$ corresponds to
\[
 \frac{1}{2\mu}G_{k,k}^{(N)}-\delta_k^{(N)}\geq 0.
\]
Finally, $\gamma^{(N)}\geq 0$ is associated with \textup{(P4)}, and $\eta_k^{(N)}\geq 0$ is associated with \textup{(P5)}. Since $t$ denotes the stepsize, the notation $\eta_k^{(N)}$ is reserved for these dual multipliers.

Define
\[
 \balpha^{(N)} :=(\alpha_0^{(N)},\ldots,\alpha_N^{(N)})^{\top},
 \qquad \bbeta^{(N)} :=(\beta_0^{(N)},\ldots,\beta_N^{(N)})^{\top}.
\]
The dual of the full pairwise SDP is analogous, but introduces one nonnegative multiplier for every ordered pair. Since our analysis starts from the consecutive-pair relaxation~\eqref{primal:reSDP}, we only state its dual:
\begin{align*}\label{Dual:reSDP}\tag{$\widehat{\mathrm D}_N$}
\min\quad & \gamma^{(N)}\Delta \\
{\rm s.t.}\quad
& \beta_0^{(N)}-\alpha_0^{(N)}-\eta_0^{(N)} -\lambda_0^{(N)}+\lambda_0^{\prime(N)}+\gamma^{(N)}=0, \\
& \beta_k^{(N)}-\alpha_k^{(N)}-\eta_k^{(N)} +\lambda_{k-1}^{(N)}-\lambda_k^{(N)} -\lambda_{k-1}^{\prime(N)}+\lambda_k^{\prime(N)}=0,~~1\leq k\leq N-1,  \tag{D1}\label{DD1} \\
& \beta_N^{(N)}-\alpha_N^{(N)}-\eta_N^{(N)} +\lambda_{N-1}^{(N)}-\lambda_{N-1}^{\prime(N)}=0,\\
& \sum_{k=0}^N\eta_k^{(N)}=1, \tag{D2}\label{DD2}\\
& \BS^{(N)} :=\sum_{k=0}^{N-1}
 \left( \lambda_k^{(N)}\BA_{k,k+1}^{(N)} +\lambda_k^{\prime(N)}\BA_{k+1,k}^{(N)}\right) 
 +\diag\left( \frac{\balpha^{(N)}}{2L}-\frac{\bbeta^{(N)}}{2\mu}\right) \succeq 0, \tag{D3}\label{DD3}\\
& \lambda_k^{(N)}\geq 0, \quad\lambda_k^{\prime(N)}\geq 0,~~ 0\leq k\leq N-1,\\
& \alpha_k^{(N)}\geq 0, \quad\beta_k^{(N)}\geq 0, \quad\eta_k^{(N)}\geq 0,~~ 0\leq k\leq N,\\
& \gamma^{(N)}\geq 0.
\end{align*}
The dual slack matrix $\BS^{(N)}$ is tridiagonal because each consecutive coefficient matrix is supported on a $2\times2$ principal block.

As shown in Lemma~\ref{Lemma:primal_dual}, Slater's condition holds under the standing assumptions. Hence $\widehat{\mathrm P}_N$ and $\widehat{\mathrm D}_N$ have the same optimal value. An explicit feasible primal--dual pair with matching objective values therefore certifies $\widehat{{\cal R}}_N(t)$. We then show that the primal solution to~\eqref{def:sdp_pep} also satisfies every nonconsecutive constraint in $\overline{\mathrm P}_N$, then  the dual matrix $\BS^{(N)}$ in~\eqref{Dual:reSDP} remains the same after adding the nonconsecutive pairs and they satisfy the KKT condition in~\eqref{def:sdp_pep}, and thus we prove $\overline{{\cal R}}_N(t) =\widehat{{\cal R}}_N(t).$

\section{Main theorems}\label{s:main}

Based on the previous discussion, we are now ready to state the main result.
Recall that
\[
{\cal R}_N(t)\leq \overline{\cal R}_N(t)\leq \widehat{\cal R}_N(t),
\]
where ${\cal R}_N(t)$ is the true worst-case performance over
$\mathcal{F}_{\mu,L}^{\PL}$, $\overline{\cal R}_N(t)$ is the value of the
full pairwise SDP~\eqref{def:sdp_pep}, and $\widehat{\cal R}_N(t)$ is the
value of the reduced SDP~\eqref{primal:reSDP}. We also define
\[
\kappa:=\frac{L}{\mu},\qquad s:=Lt-1.
\]

\begin{theorem}[\bf Main theorem]
\label{thm:main}
Let $0<t<2/L$, and consider $N$ steps of gradient descent
\[
\bx_{k+1}=\bx_k-t\nabla f(\bx_k),\qquad 0\leq k\leq N-1,
\]
with constant stepsize $t$ and iterates $(\bx_0,\bx_1,\ldots,\bx_{N})$. Then the consecutive-pair and full pairwise SDP values satisfy the following statements.

\begin{enumerate}
\item  \text{\bf Regimes I and II: $0<t\leq \sqrt{3}/L$.}
For all $N\geq 1$,
\[
\widehat{\cal R}_N(t)=p(t)^N,
\]
where
\begin{equation}\label{def:p}
p(t)= \begin{cases}
\displaystyle
\frac{s^2+2(\kappa-1)}{2\kappa+1}
\left( \frac{ \sqrt{2(\kappa+1)s^2+2(\kappa-1)(2\kappa+1)}+|s|}{\sqrt{2(\kappa+1)s^2+2(\kappa-1)(2\kappa+1)}-|s|}\right),
& 0<t\leq \dfrac{1}{L}, \\
\displaystyle
1+\frac{(s^2-1)(s+3)}{s+2\kappa+1}, & \dfrac{1}{L}<t\leq \dfrac{\sqrt{3}}{L}.
\end{cases}
\end{equation}
Moreover, in this range the constructed optimal primal solution
of the reduced SDP~\eqref{primal:reSDP} is also feasible for the
full pairwise SDP~\eqref{def:sdp_pep}. The corresponding reduced dual certificate extends to the full dual by assigning zero multipliers to the nonconsecutive constraints. Hence
\[
\overline{\cal R}_N(t) = \widehat{\cal R}_N(t) = p(t)^N.
\]
Therefore, for every $f\in\mathcal{F}_{\mu,L}^{{\rm \PL}}$ and $0 <Lt \leq \sqrt{3}$, it holds that
\[
f(\bx_{N})-f^{\star} \leq p(t)^N\left(f(\bx_0)-f^{\star}\right).
\]

\item \text{\bf Regime III: $\sqrt{3}/L<t < 2/L$.}
Define
\[
p(t):=1+\frac{(s^2-1)(s+3)}{s+2\kappa+1},~~~\frac{\sqrt{3}}{L} < t < \frac{2}{L}
\]
which extends the expression from Regime~II, and let
\begin{equation}\label{def:c_N0}
 c^{(N)} := \frac{\kappa s^{2N}}{1+(\kappa-1)s^{2N}},
 \qquad N_0 := \left\lfloor \frac{\log\!\left(\frac{2}{s+3}\right)}{2\log s} \right\rfloor.
\end{equation}
The following two phases arise.

\begin{enumerate}
\item[\textup{(a)}]
\emph{Finite-horizon exactness for $1\leq N\leq N_0$.}
The constructed optimal primal solution of the consecutive-pair
SDP also satisfies every nonconsecutive interpolation constraint
of the full pairwise SDP. Thus it holds that
\[
\overline{\cal R}_N(t) = \widehat{\cal R}_N(t) = c^{(N)}, \qquad 1\leq N\leq N_0.
\]

\item[\textup{(b)}]
\emph{Geometric tail and asymptotic exactness for $N>N_0$.}
The consecutive-pair SDP has the exact optimal value
\[
\widehat{\cal R}_N(t) = c^{(N_0)}p(t)^{N-N_0}.
\]
Moreover, the Toeplitz-envelope primal construction from
Regime~II extends to this stepsize range and is feasible for
the full pairwise SDP, with objective value $\Delta p(t)^N$.
Hence
\[
p(t)^N \leq \overline{\cal R}_N(t) \leq \widehat{\cal R}_N(t)
= c^{(N_0)}p(t)^{N-N_0}, \qquad N>N_0.
\]
In particular, for each fixed $t\in(\sqrt{3}/L,2/L)$,
\[
\lim_{N\to\infty}
\overline{\cal R}_N(t)^{1/N} = \lim_{N\to\infty} \widehat{\cal R}_N(t)^{1/N} = p(t).
\]
Thus $p(t)$ is the exact asymptotic geometric factor of both
SDP relaxations. No finite-horizon equality between their
optimal values is asserted for $N>N_0$.
\end{enumerate}

In conclusion, for every $f\in{\cal F}_{\mu,L}^{\rm \PL}$ and $\sqrt{3} < Lt < 2$, we have
\[
f(\bx_N)-f^{\star}
\leq
\begin{cases}
c^{(N)}\left(f(\bx_0)-f^{\star}\right),
& 1\leq N\leq N_0,\\[1mm]
c^{(N_0)}p(t)^{N-N_0}
\left(f(\bx_0)-f^{\star}\right),
& N>N_0.
\end{cases}
\]
\end{enumerate}
The exactness statements above concern the SDP relaxations~\eqref{primal:reSDP};
they do not assert tightness over the full function class
${\cal F}_{\mu,L}^{\PL}$. In particular, when $N=1$, the bound reduces to the one-step result
in~\cite{AKZ23}.
\end{theorem}

\begin{figure}[h!]
\centering
\includegraphics[width=120mm]{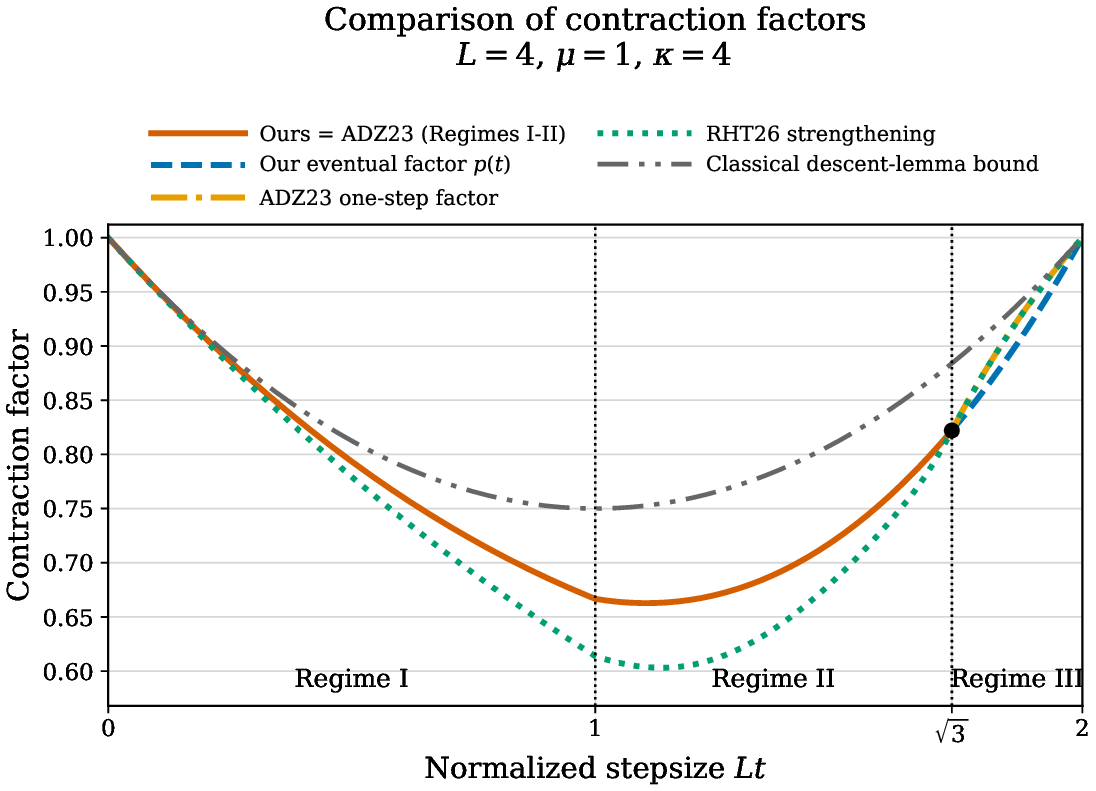}
\caption[Comparison of classical and SDP-based contraction factors]
{Comparison of contraction factors for constant-stepsize gradient
descent with $L=4$, $\mu=1$, and $\kappa=L/\mu=4$.
The horizontal axis is the normalized stepsize $Lt$.
The eventual factor $p(t)$ in Theorem~\ref{thm:main} coincides
with the one-step SDP factor of~\protect\cite{AKZ23} for
$0<Lt\leq\sqrt{3}$ (solid red curve), whereas for
$\sqrt{3}<Lt<2$ it is strictly smaller (blue dashed versus
orange dash-dotted curves).
The green dotted curve shows the numerically computed one-step
bound obtained using the one-point strengthening
of~\protect\cite{RHT26}.
The gray long-dash--double-dot curve represents the classical
descent-lemma bound
$p_{\rm van}(t)=1-\mu t(2-Lt)$;
see~\protect\cite{KNS16}.
The green dot curve corresponds to a strengthened SDP~\cite{RHT26}, distinct from the relaxations whose exact values are characterized here. Smaller factors indicate stronger convergence guarantees; SDP exactness does not imply tightness over the smooth P\L{} function class because of the lack of exact interpolation condition.}
\label{fig:compare}
\end{figure}

Theorem~\ref{thm:main} shows that the convergence behavior depends qualitatively on the stepsize $t$, with three regimes:
\[
0<t\leq \frac{1}{L},\qquad
\frac{1}{L}<t\leq \frac{\sqrt{3}}{L},\qquad
\frac{\sqrt{3}}{L}<t<\frac{2}{L}.
\]
We visualize the corresponding convergence factors in
Figure~\ref{fig:compare}, where $p(t)$ denotes the eventual linear convergence rate.

\paragraph{Comparisons with existing works}
As a classical non-PEP baseline, combining the descent lemma
with the P\L{} inequality gives the $N$-step contraction bound
$p_{\rm van}(t)^N$, where
$p_{\rm van}(t)=1-\mu t(2-Lt)$ for $0<t<2/L$;
see~\cite{KNS16} for the standard choice $t=1/L$.
Our SDP analysis sharpens this baseline through the
finite-horizon guarantees in Theorem~\ref{thm:main},
with the improvement in contraction factors illustrated
in Figure~\ref{fig:compare}.

Our result extends the one-step analysis of~\cite{AKZ23} to the multi-step setting. In particular, when $N=1$, our formulas reduce to the same one-step expressions as in~\cite{AKZ23}. The multi-step analysis leads to two distinct phenomena.

\begin{figure}[h!]
  \centering
 \includegraphics[width=120mm]{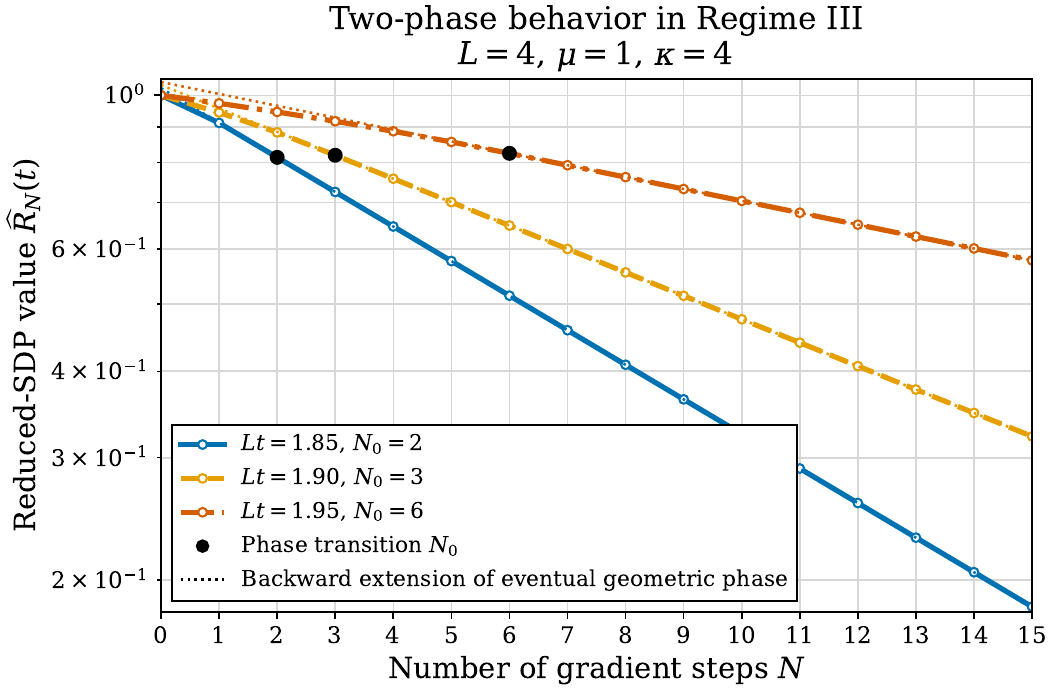}
 \caption{Two phases for regime $3$, $x$-axis: number of timestep $N$, $y$-axis: objective value of~\eqref{primal:reSDP}. The convergence is not linear at the first phase but becomes linear at the second phase. The linear rate at second phase is the same as the red line in Figure~\ref{fig:compare}. For this example plot, we set $L=4$ and $\mu=1$.}
 \label{fig:phase_3}

\end{figure}

\begin{enumerate}
\item
For the first two regimes, namely $0<t\leq 1/L$ and $1/L<t\leq \sqrt{3}/L$, we recover exactly the same convergence factor as in~\cite{AKZ23}, but now for every horizon $N$. Hence the one-step bound from~\cite{AKZ23} \emph{tensorizes} in the
multi-step setting:
\[
\widehat{\cal R}_N(t)=\overline{\cal R}_N(t)=p(t)^N.
\]
In particular, this proves that the corresponding rate is tight at the
level of the full pairwise SDP.

\item
For the third regime, namely $\sqrt{3}/L<t<2/L$, the
consecutive-pair SDP exhibits two horizon-dependent phases.
During the transient phase $1\leq N\leq N_0$, where $N_0$
is defined in~\eqref{def:c_N0}, the optimal value has a
nongeometric dependence on the horizon, and the constructed
primal optimizer is also feasible for the full pairwise SDP.
It holds that
\[
\overline{\cal R}_N(t) =\widehat{\cal R}_N(t)
=c^{(N)}, \qquad 1\leq N\leq N_0.
\]
For $N>N_0$, the consecutive-pair SDP value follows a
geometric decay pattern:
\[
\widehat{\cal R}_N(t) =c^{(N_0)}p(t)^{N-N_0}
\]
where $c^{(N)}(t)$ is defined in~\eqref{def:c_N0}. The particular mixed primal optimizer used in this phase violates some nonconsecutive constraints of the full pairwise SDP~\eqref{def:sdp_pep}. 

Nevertheless, a separate full-SDP feasible construction gives
\[
p(t)^N \leq\overline{\cal R}_N(t) \leq\widehat{\cal R}_N(t) =c^{(N_0)}p(t)^{N-N_0}.
\]
Thus both SDP relaxations have the same asymptotic geometric factor $p(t)$, which is strictly smaller than the factor obtained by recursively applying the one-step bound
of~\cite{AKZ23} to the same stepsize regime. Figure~\ref{fig:phase_3} illustrates the transition from the nongeometric transient to the geometric tail of the consecutive-pair SDP value.
\end{enumerate}

Specializing the proximal-gradient bound of~\cite{KJH26} to a vanishing nonsmooth term gives
\[
p_{\rm KJH}(t) = \begin{cases}
\displaystyle
\frac{\kappa(Lt+1)^2}
{\kappa(Lt+1)^2+Lt(Lt+2)}, &0<Lt\leq\sqrt{3}, \\ 
\displaystyle
\frac{\kappa(Lt-1)^2}
{\kappa(Lt-1)^2+Lt(2-Lt)}, & \sqrt{3}<Lt<2.
\end{cases}
\]
For $1<Lt<\sqrt{3}$, direct calculation gives
\[
p_{\rm KJH}(t)-  p_{\rm ADZ}(t) = \frac{ L^2t^2(\kappa+1)(Lt+2)(3-L^2t^2) }{ (Lt+2\kappa) \bigl[\kappa(Lt+1)^2+Lt(Lt+2)\bigr] } >0,
\]
where $  p_{\rm ADZ}(t)$ denotes the direct smooth-P\L{} one-step factor
from~\cite{AKZ23}. On the other hand, for $\sqrt{3}\leq Lt<2$, it holds that
\[
p_{\rm KJH}(t) = \frac{\kappa(Lt-1)^2} {1+(\kappa-1)(Lt-1)^2} =  p_{\rm ADZ}(t).
\]
Thus on $1<Lt<\sqrt{3}$, the general composite certificate yields a valid but strictly weaker smooth-gradient-descent bound, whereas the two bounds coincide in the large-stepsize regime.

Another relevant work is by Rubbens, Hendrickx and Taylor~\cite{RHT26} who derive a tighter inequality for $\mathcal{F}^{\PL}_{\mu,L}$ by one-point strengthening technique. They propose to replace the smoothness interpolation inequality~\eqref{eq:PLsmooth} with
\begin{equation}
\label{eq:PLsmooth-oss}
\begin{aligned}
& f_i-f_j-\lag\bg_j,\bx_i-\bx_j\rag
 -\frac{1}{2L}\|\bg_i-\bg_j\|^2
 +\frac{L}{4}
 \left\|\bx_i-\bx_j-\frac{1}{L}(\bg_i-\bg_j)\right\|^2  \\
 & \geq \max_{0 \leq \alpha \leq \frac{2\mu}{2L+\mu}}  \frac{\alpha}{(1-\alpha)(2\mu-(L+\mu)\alpha)} \\
 & \qquad \cdot \left((1-\alpha)^2(L+\mu)\left(f_i-f^*-\frac{\|\bg_i\|^2}{2L}\right)-(L-\mu)(f_j-f^*+\frac{\|\bg_j\|^2}{2L}) \right).
\end{aligned}
\end{equation}
It is straightforward to see that inequality~\eqref{eq:PLsmooth-oss} is stronger than inequality~\eqref{eq:PLsmooth}. The PEP bound with one-point strengthening technique indeed yield a tighter bound for one step formulation when $t<\sqrt{3}/L$~(see green line in Figure~\ref{fig:compare}). However, due to the complicated form of the one-point strengthening inequality, the convergence rate is computed via numerics. We evaluate the strengthened SDP numerically while no closed-form value is derived here.

\paragraph{Summary of the theorems for each regime}
The proof is divided into several parts based on each regime of the stepsize.  By Theorem~\ref{thm:primal_dual_equal}, and~\ref{thm:primal_dual_equal_case2} for $0< Lt \leq 1$ and $1<Lt\leq \sqrt{3}$ respectively, we have
\[
\eqref{primal:reSDP}=\Delta\widehat{\cal R}_N(t)=\Delta p(t)^N.
\]
Moreover, Theorem~\ref{thm:primal_dual_equal} and~\ref{thm:primal_dual_equal_case2} show that,  the optimal primal and dual solutions of the reduced SDP also satisfy all constraints of the full SDP~\eqref{def:sdp_pep} in these two regimes. Therefore, we have
\[
{\cal R}_N(t)\leq \overline{\cal R}_N(t)=\widehat{\cal R}_N(t)=p(t)^N,
\]

For $\sqrt{3}<Lt<2$, Theorem~\ref{thm:primal_dual_equal_case31} and~\ref{thm:primal_dual_equal_case32} imply that
\[
\widehat{\cal R}_N(t)=
\begin{cases}
c^{(N)}, & N\leq N_0,\\[1mm]
c^{(N_0)}\,p(t)^{\,N-N_0}, & N>N_0.
\end{cases}
\]
For transient phase $1\leq N\leq N_0$, Theorem~\ref{thm:primal_dual_equal_case31} additionally asserts
\[
\overline{\cal R}_N(t)=\widehat{\cal R}_N(t)=c^{(N)}.
\]
For geometric phase $N>N_0$, Proposition~\ref{prop:r3c2:toeplitz-full-feasible} shows $p(t)^N \leq \overline{\cal R}_N(t).$
Since $$\mathcal{R}_N(t)  \leq \overline{\cal R}_N(t)\leq \widehat{\cal R}_N(t),$$ we have the corresponding upper bound for $f(\bx_{N})-f^{\star}$ follows
immediately.

\paragraph{Brief discussion on the constructions of primal and dual certificates}
For each regime, we explicitly construct primal and dual
feasible solutions for the reduced SDP~\eqref{primal:reSDP}
and its dual~\eqref{Dual:reSDP}, with matching objective
values. The duality then certifies their optimality. The detailed constructions of the primal/dual certificates for each regime are given in the next subsections. 
In the first two regimes, i.e., $0\leq Lt\leq 1$ and $1\leq Lt\leq \sqrt{3}$, these optimal solutions also satisfy
all constraints of the full SDP~\eqref{def:sdp_pep}, which yields
\[
\overline{\cal R}_N(t)=\widehat{\cal R}_N(t).
\]
In Regime~III ($\sqrt{3}\leq Lt\leq 2$), the constructed primal optimizer remains feasible for the full pairwise SDP during the transient phase $1\leq N\leq N_0$. For $N>N_0$, the particular mixed primal optimizer in~\eqref{primal:reSDP} violates some nonconsecutive constraints in~\eqref{def:sdp_pep}. The transition between the two constructions is determined by the positive-semidefiniteness threshold of the transient dual slack matrix $\BS^{(N)}$.

\begin{figure}[h!]
\centering
\includegraphics[width=160mm]{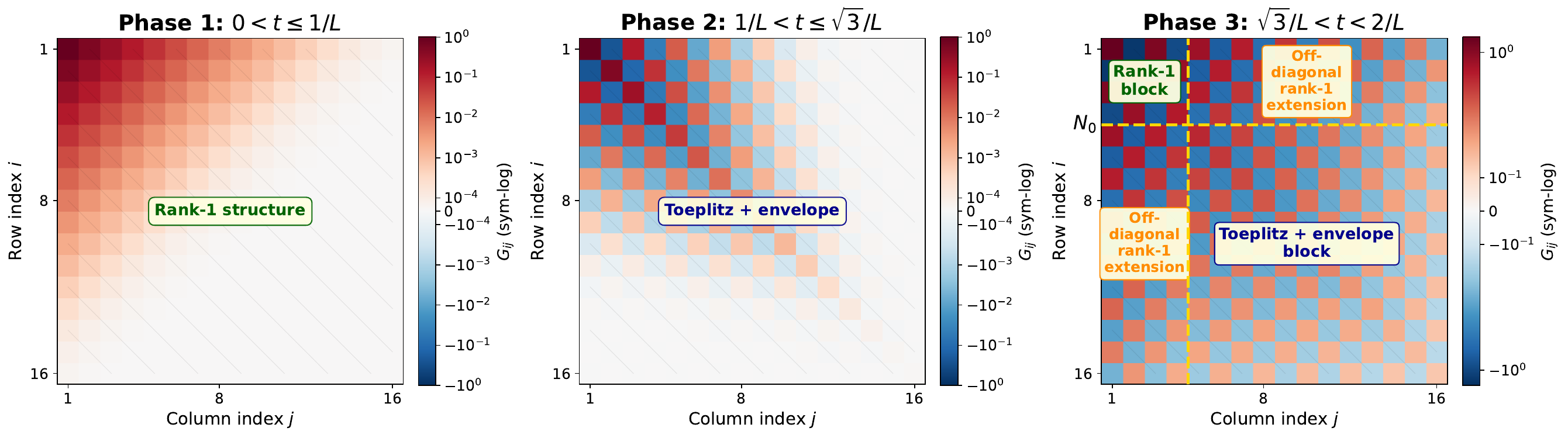}
\caption{Structure of the constructed optimal Gram matrix
$\BG$ for the consecutive-pair SDP~\eqref{primal:reSDP} in Regime I-III.
From left to right: the rank-one construction in Regime~I;
the diagonally scaled Toeplitz construction in Regime~II;
and the mixed rank-one-prefix and Toeplitz-tail construction
in the geometric phase of Regime~III, where $N>N_0$.}
\label{fig:gstructure}
\end{figure}

The structure of these primal and dual solutions is computer-assisted:
we first infer the pattern from numerical SDP solutions and then verify it
rigorously through explicit constructions. The solution structure of $\BG$ is illustrated in Figure~\ref{fig:gstructure}.

The constructed Gram matrices exhibit three distinct structures. In Regime~I, $\BG$ is the rank-one matrix
in~\eqref{psol:reSDP-r1}: the sampled gradients are collinear, and their squared norms decay geometrically with ratio $p(t)$. In Regime~II, the construction~\eqref{psol:reSDP-r2} takes
the form $\BG=2\mu\Delta\BP^{(N)}\BPhi^{(N)}\BP^{(N)}$,
where the Toeplitz matrix $\BPhi^{(N)}$ describes gradient
correlations depending only on the index difference, while
the diagonal factors $\BP^{(N)}$ impose a geometric envelope.
The resulting Gram matrix is positive definite.
In Regime~III, the transient construction
in~\eqref{psol_r3c1:reSDP} remains rank one, with consecutive
gradients pointing in opposite directions.
For $N>N_0$, the mixed construction
in~\eqref{psol_r3c2:reSDP} joins this rank-one prefix to a
diagonally scaled Toeplitz tail. Its off-diagonal blocks
preserve the relation that each prefix gradient is a scalar
multiple of the gradient at the transition index $N_0$.
Figure~\ref{fig:gstructure} illustrates the Regime-I,
Regime-II, and mixed Regime-III constructions.

\subsection{Construction of the exact optimal solution for the regime $t\leq 1/L$}
We now construct the certificates underlying the main result.
For each regime, we exhibit primal and dual feasible points with
matching objective values, which certify optimality by weak duality.
Whenever the primal point also satisfies the omitted nonconsecutive
constraints, it is optimal for the full pairwise SDP as well.
The distinction between these SDP certificates and interpolation
by an actual smooth PL function will be important below.
We begin with $0<t\leq1/L$.

For $t \leq 1/L$, the solution for the relaxed primal problem~\eqref{primal:reSDP}
and~\eqref{Dual:reSDP} is characterized by the following theorem.

\begin{theorem}[\bf Regime I]\label{thm:primal_dual_equal}
Let
$$
0<t \leq 1/L,
$$ and let $N \geq 1$ be an integer, then the following statements hold.
\begin{enumerate}
    \item  \label{thm:primal_dual_equal_case_r1} An optimal primal solution of the consecutive-pair SDP~\eqref{primal:reSDP} is given by
\begin{equation}\label{psol:reSDP-r1}
\begin{aligned}
\delta_k &= \Delta p^{k},~~0 \leq k\leq N, \\
\ell &= \Delta p^{N},~
\BG = 2\mu \Delta 
\begin{bmatrix}
1 \\
\vdots \\
p^{N/2}
\end{bmatrix}
\begin{bmatrix}
1 \\
\vdots \\
p^{N/2}
\end{bmatrix}^{\top}.
\end{aligned}
\end{equation}
A corresponding dual solution to~\eqref{Dual:reSDP} is given by
\begin{equation}\label{dsol:reSDP}
\begin{aligned}
&\alpha_k  = 0,~~0 \leq k\leq N,   \\
&\beta_0 =p^{N}(rp^{-1}-1),~~\beta_{N} = 1-r, \\
&\beta_k  = rp^{N-k}(p^{-1}-1),~~ ~1 \leq k\leq N-1, \\
&\lambda_k  = rp^{N-k-1}, \lambda'_{k}=0,~~ ~ 0 \leq k\leq N-1,\\
& \gamma = p^{N}, \\
& \eta_{N} = 1,~~\eta_k  = 0,~~0\leq k\leq N-1,
\end{aligned}
\end{equation}
where 
\begin{equation}\label{eq:def_rho}
r:= \frac{2 \kappa (s^2+2\kappa-2) }{s^2+(2\kappa+1)(s^2+2\kappa-2) - \sqrt{2(\kappa+1)s^2+2(\kappa-1)(2\kappa+1)}|s|}
\end{equation}
and $p$ is defined in~\eqref{def:p}. 
    \item  \label{thm:primal_dual_equal_case_r2} We have the solution is also an optimal solution of~\eqref{def:sdp_pep}.
\end{enumerate}
\end{theorem}
\noindent\textbf{Proof insights.}
The primal construction uses a rank-one Gram matrix $\BG$ whose
diagonal entries decay geometrically with ratio $p$. The objective
gaps $\delta_k$ follow the same geometric pattern, and imposing
equality in~\eqref{PP1} determines the contraction factor $p$.

For the dual construction, complementary slackness suggests setting
$\alpha_k=0$ and $\lambda'_k=0$. The interior multipliers
$\{\beta_k\}_{k=1}^{N-1}$ and the multipliers
$\{\lambda_k\}_{k=0}^{N-1}$ each form a geometric sequence with
ratio $1/p$. This structure, together with the dual balance
equations~\eqref{DD1}, yields explicit expressions for the
multipliers.

Finally, we verify the positive semidefiniteness condition~\eqref{DD3}
by expressing $\BS$ as a nonnegative weighted sum of rank-one
positive semidefinite matrices; see Lemma~\ref{R1:L3}.\\
\begin{proof}[\bf Proof for Theorem~\ref{thm:primal_dual_equal}]

The proposed primal and dual points have matching objective values, $\ell=\gamma\Delta=p^N\Delta$. By weak duality in Lemma~\ref{Lemma:primal_dual}, it remains only to verify primal and dual feasibility. \\ 
\\
\noindent\textbf{Nonnegativity.}
We first verify that all dual multipliers are nonnegative.
By Lemma~\ref{R1:L1}, we have $0<p<r<1$.
The multipliers $\alpha_k$, $\eta_k$, and $\lambda'_k$
are nonnegative by definition. Moreover,
\[
\lambda_k=rp^{N-k-1}>0,
\qquad 0\leq k\leq N-1,
\qquad
\gamma=p^N>0.
\]
For the multipliers $\beta_k$, we obtain
\begin{align*}
\beta_0
&=p^N(rp^{-1}-1)
 =p^{N-1}(r-p)>0,\\
\beta_k
&=rp^{N-k}(p^{-1}-1)
 =rp^{N-k-1}(1-p)>0,
 \qquad 1\leq k\leq N-1,\\
\beta_N
&=1-r>0.
\end{align*}
Thus all dual multipliers are nonnegative.\\
\\
\textbf{Primal Feasibility.}  We verify primal constraints~\eqref{PP1}--\eqref{PP6} in order. \\
\\
\eqref{PP1},~\eqref{PP2}: By Lemma~\ref{R1:L4} , the constraints all hold. In particular, the forward consecutive pairs $(k,k+1)$ for $0 \leq k \leq N-1$ hold with equality. \\
\\
\eqref{PP3}: By definition, we have 
$$
\frac{G_{k,k}}{2\mu}  = p^{k}\Delta = \delta_k > \frac{G_{k,k}}{2L}.\\
$$
\\
~\eqref{PP4}: By definition, we have $\delta_0$ = $\Delta$.\\
\\           
\eqref{PP5}: Since $p<1$, we have
\[
\ell = \delta_{N} < \delta_{N-1} <\cdots<\delta_0
\]
since $\delta_{k+1} = p\delta_{k}, 0\leq k\leq N-1$ by definition.\\
\\
\eqref{PP6}: Notice we have $\BG$ is a rank-one matrix, i.e. $\BG =\bv\bv^{\top} \succeq 0$ where $\bv\in\RR^{N+1}$ where
\[
v_k =  \sqrt{2\mu \Delta p^k},~0\leq k\leq N.
\]
\textbf{Dual Feasibility.} We proceed to verify dual constraints~\eqref{DD1}--\eqref{DD3} in order.\\
\\
\eqref{DD1}: The equality constraints are easily verified by directly plugging in the expressions.\\
\\
\eqref{DD2}: The constraint holds trivially since $\eta_{N}= 1$ and $\eta_{k}=0~(k=0,\dots N-1)$.\\
\\
\eqref{DD3}: The constraint is directly verified by~Lemma~\ref{R1:L3}.\\
\\
\textbf{Full-SDP Feasibility.} We proceed with item 2 in Theorem~\ref{thm:primal_dual_equal}. Since~\eqref{primal:reSDP} is a relaxation of~\eqref{def:sdp_pep}, it remains to verify $\BG$ and $\left\{\delta_k\right\}$ are feasible for non-consecutive pairs for the first constraint in~\eqref{def:sdp_pep}. This condition is directly verified by Lemma~\ref{R1:L4}. Extend the reduced dual certificate by assigning zero multipliers
to the omitted nonconsecutive constraints. The resulting full-SDP
primal and dual points remain feasible and have matching objective
values. Weak duality therefore establishes their optimality for
the full pairwise SDP.
\end{proof}
In this regime, the optimal Gram matrix $\BG$ has rank one, so the
sampled gradients are collinear. This suggests seeking a
one-dimensional function that interpolates the optimal data.
However, the following proposition shows that no $\mu$-PL,
$L$-smooth function on $\mathbb{R}$ can generate the data
in~\eqref{psol:reSDP-r1} through gradient descent with stepsize $t$.

The key observation is that, after normalizing $f^\star=0$ and
choosing the orientation so that the sampled gradients are positive,
smoothness ensures that $f'$ remains positive between consecutive
iterates. The PL inequality therefore gives
$f'(x)\geq\sqrt{2\mu f(x)}$ on these intervals. Integrating this
inequality forces the objective gaps to contract by a factor
strictly smaller than $p$, contradicting the prescribed
geometric decay.
\begin{proposition}\label{PROP:EX}
Let $0 < t \leq 1/L$ and $N \geq 1$ be an integer. Let $\delta,\BG$ be the optimal solution given by~\eqref{psol:reSDP-r1}. We have there is no $\mu$-PL and $L$-smooth function $f$ in 1D such that
\begin{align*}
& f(x_k)-f^* = \delta_k, ~~f'(x_k) = g_k ~~0 \leq k \leq N, \\ & x_{k+1} = x_k - tf'(x_k),~~0 \leq k \leq N-1,
\end{align*}
where
\[
\BG = \bg\bg^{\top}, ~~~ \bg = [g_0,g_1,\cdots,g_N]^{\top}.
\]
\end{proposition}
\begin{proof}[\bf Proof for Proposition~\ref{PROP:EX}]
We prove by contradiction. Suppose such function $f$ exists, since $f$ is $\mu$-PL, we normalize the minimum by assuming $f^{\star} = 0$  and $g_i>0$ without loss of generality, hence $f \geq 0$, then we have
$$
f'(x)^2  \geq 2\mu f(x).
$$
\textbf{Positivity of derivative.} Take any interval $[x_{i+1},x_i]$, since $f$ is $L$-smooth we have for any $x \in  (x_{i+1},x_i)$ 
$$
f'(x) \geq f'(x_i) - L (x_i-x) > g_i - L t g_i  \geq 0
$$
since $ 0 < t \leq 1/L$. Hence we have $f'(x)$ must be positive on $[x_{i+1},x_i]$, then $f'(x) \geq \sqrt{2\mu f(x)}$ and $ f(x_i) \geq f(x) \geq f(x_{i+1}) >0$. \\
\\
\textbf{Integrated PL inequality.} We define $u(x) = \sqrt{f(x)}$, then we have
$$
u'(x) = \frac{f'(x)}{2\sqrt{f(x)}} \geq \frac{\sqrt{2\mu f(x)}}{2\sqrt{f(x)}} = \sqrt{\frac{\mu}{2}}.
$$
Hence
$$
u(x_{i}) - u(x_{i+1})=\int^{x_{i}}_{x_{i+1}} u'(x) dx \geq \sqrt{\frac{\mu}{2}} (x_i-x_{i+1}) =  \sqrt{\frac{\mu}{2}}t g_i = \mu t\sqrt{f_i}.
$$
We necessarily have
$$
(1-\sqrt{p}) \sqrt{f_i} \geq \mu t \sqrt{f_i} \Longleftrightarrow \sqrt{p}  \leq 1- \mu t.
$$
\textbf{Contradiction.} Recall $\sqrt{p}$ is the positive root of the quadratic equation $h(x)$~\eqref{eq:psolqua}, hence we evaluate at 
$$
h(1-\mu t):= (2\kappa+1)(1-\mu t)^2 + 2s(1-\mu t)+(2-s^2-2\kappa)= (\mu^2-L^2)t^2 < 0.
$$
This implies $\sqrt{p}>1-\mu t$, which contradicts the necessary inequality $\sqrt{p} \leq 1-\mu t$.
\end{proof}
This proposition concerns realization of the prescribed SDP data by a one-dimensional gradient-descent trajectory; it does not invalidate the SDP optimality certificate. Nor does the  rank one property of the sampled Gram matrix alone reduce interpolation in higher
dimensions to a one-dimensional problem, since the gradient between sampled points need not remain collinear. We next turn to Regime~II, where the certificate has a positive definite Gram matrix and requires a different construction.
\subsection{Constructions of the exact optimal solution for the regime $1 < Lt \leq \sqrt{3}$}
\begin{theorem}[\bf Regime II]\label{thm:primal_dual_equal_case2}
Let 
\[
1/L < t \leq \sqrt{3}/L,
\]
and let $N \geq 1$ be an integer, then the following statements hold.
\begin{enumerate}
    \item  \label{thm:primal_dual_equal_case2_r1} An optimal primal solution to~\eqref{primal:reSDP} is given by
\begin{equation}\label{psol:reSDP-r2}
\begin{aligned}
\delta_k &= p^{k}\Delta,~\forall~0\leq k\leq N,~\ell = p^{N}\Delta, \\
\BG &= 2\mu\Delta \BP^{(N)} \BPhi^{(N)} \BP^{(N)},
\end{aligned}
\end{equation}
where
\begin{equation}\label{def:BPBPHI}
\begin{aligned}
\BP^{(N)}
=\diag\left(1,p^{1/2},\ldots,p^{N/2}\right),
\qquad
\BPhi^{(N)}
=\left(\phi_{|i-j|}\right)_{i,j=0}^{N},
\end{aligned}
\end{equation}
and $\phi_k$'s are given by the following recursive relation
\begin{equation}\label{R2:eq_recur}
\phi_0=1,~~\phi_1=\frac{1}{2\sqrt{p}}(1-L^2t^2+p)=\frac{\tau}{\sqrt{p}},~~\phi_{k+2} = \frac{s}{s+2}\phi_k + \frac{p-s(s+2)}{\sqrt{p}(s+2)}\phi_{k+1},~k\geq 0
\end{equation}
with $p$ defined by~\eqref{def:p} and
\[
\tau=\frac{1}{2}(1-L^2t^2+p) = \frac{2-L^2t^2-\mu t}{2+\mu t}.
\]
The corresponding dual variables to~\eqref{Dual:reSDP} given by
\begin{align*}
&\lambda_{k} = p^{N-k-1} \frac{s+2}{2+\mu t},~~ \lambda_k' = p^{N-k-1} \frac{s}{2+\mu t},~~0\leq k\leq N-1, \\
&\alpha_{k} = 0,~~0\leq k\leq N,\\
&\eta_{N} =1,~~\eta_k = 0,~~0\leq k\leq N-1,\\
&\gamma = p^N,\\
&\beta_k = \lambda'_{k-1}-\lambda'_k - (\lambda_{k-1}-\lambda_k),~1\leq k\leq N-1,\\
&\beta_0 = \lambda_0-\lambda_0'-p^N,~~\beta_{N} = 1 + \lambda_{N-1}' - \lambda_{N-1}.
\end{align*}
\item \label{thm:primal_dual_equal_case2_r2} We have the solution is also an optimal solution of~\eqref{def:sdp_pep}.
\end{enumerate}
\end{theorem}
\noindent\textbf{Proof insights.}
The main challenge in this regime is to construct a full-rank Gram
matrix $\BG$, since the consecutive-pair constraints leave its entries
outside the tridiagonal band undetermined.

We begin with the dual construction. The full-rank ansatz for $\BG$
suggests imposing $\BS=0$, as required by complementary slackness
at optimality. We also choose the multipliers
$\{\lambda_k\}_{k=0}^{N-1}$ and
$\{\lambda'_k\}_{k=0}^{N-1}$ to form geometric sequences with
ratio $1/p$. Combining this structure with $\BS=0$ and the dual
balance equations~\eqref{DD1} yields explicit expressions for the
dual multipliers.

For the primal construction, we choose the objective gaps $\delta_k$
and the diagonal entries of $\BG$ to decay geometrically with
ratio $p$, and impose equality in~\eqref{PP1} and~\eqref{PP2}.
These conditions determine the consecutive gradient inner products
but leave the remaining off-diagonal entries unspecified.
To determine these entries, we additionally impose
\[
\Tr(\BA_{i,j}\BG)=\delta_i-\delta_j,
\qquad 0\leq j<i\leq N.
\]
These equalities lead to a recurrence for the remaining entries
and yield a diagonally scaled Toeplitz matrix. The resulting
structure allows us to establish $\BG\succ0$; see
Lemma~\ref{R2:L3}.

\begin{proof}[\bf Proof for Theorem~\ref{thm:primal_dual_equal_case2}]
We start with verifying the primal and dual variables.
The proposed primal and dual points have matching objective values, $\ell=\gamma\Delta=\Delta p^N$. By weak duality in Lemma~\ref{Lemma:primal_dual}, it remains only to verify all primal and dual constraints. \\
\\
\noindent\textbf{Nonnegativity.}
We first verify that all dual multipliers are nonnegative.
By Lemma~\ref{R2:L1}, we have $0<p<1$.
Moreover, $s=Lt-1>0$ in this regime.
The multipliers $\alpha_k$ and $\eta_k$ are nonnegative by
definition, while
\begin{align*}
\lambda_k
&=\frac{s+2}{2+\mu t}\,p^{N-k-1}>0,
~~0\leq k\leq N-1,\\
\lambda'_k
&=\frac{s}{2+\mu t}\,p^{N-k-1}>0,
~~0\leq k\leq N-1,\\
\gamma&=p^N>0.
\end{align*}
For the multipliers $\beta_k$, we obtain
\begin{align*}
\beta_0
&=p^{N-1}\left(\frac{2}{2+\mu t}-p\right)
 =\frac{\mu t(3-L^2t^2)}{2+\mu t}\,p^{N-1}
 \geq0,\\
\beta_k
&=\frac{2(1-p)}{2+\mu t}\,p^{N-k-1}>0,
 \qquad 1\leq k\leq N-1,\\
\beta_N
&=\frac{\mu t}{2+\mu t}>0.
\end{align*}
The inequality for $\beta_0$ follows from $Lt\leq\sqrt{3}$,
with equality at $Lt=\sqrt{3}$.
Thus all dual multipliers are nonnegative.\\
\\
\textbf{Primal Feasibility.}  We verify primal constraints~\eqref{PP1}--\eqref{PP6} in order.\\
\\
\eqref{PP1},\eqref{PP2}: These constraints hold by Lemma~\ref{R2:L2} statement~\ref{R2:L2S1} and~\ref{R2:L2S2} respectively. \\
\\
\eqref{PP3}: We have $G_{k,k} =2\mu \delta_k \leq 2L \delta_k$ by definition.\\
\\
\eqref{PP4}: By definition we have $\delta_{0} = \Delta$. \\
\\
\eqref{PP5}:We have
$$
\ell = \delta_{N} <\delta_{N-1} < \cdots<\delta_0
$$
since $p<1$ and $\delta_{k+1} = p\delta_{k}, 0\leq k \leq N-1$ by definition. \\
\\
\eqref{PP6}: We have $\BG \succeq 0$ by Lemma~\ref{R2:L3}.\\
\\
\textbf{Dual Feasibility.}  We proceed to verify the dual constraints~\eqref{DD1}--\eqref{DD3} in order. \\
\\
\eqref{DD1}: Substituting the proposed multipliers verifies the constraints. \\
\\
\eqref{DD2}: The constraint holds trivially since $\eta_{N}= 1$ and $\eta_{k}=0~(k=0,\dots N-1)$.\\
\\
\eqref{DD3}: The constraint holds since $\BS=0$ by Lemma~\ref{R2:L4}.\\
\\
\textbf{Full-SDP Feasibility.} We proceed with item 2 in Theorem~\ref{thm:primal_dual_equal_case2}. Since~\eqref{primal:reSDP} is a relaxation of~\eqref{def:sdp_pep}, it remains to verify $\BG$ and $\left\{\delta_k\right\}$ are feasible for non-consecutive pairs for the first constraint in~\eqref{def:sdp_pep}. This condition is directly verified by (b) \& (c) in~Lemma~\ref{R2:L2}. Extend the reduced dual certificate by assigning zero multipliers
to the omitted nonconsecutive constraints. The resulting full-SDP
primal and dual points remain feasible and have matching objective
values. Weak duality therefore establishes their optimality for
the full pairwise SDP.
\end{proof}
The restriction $Lt\leq\sqrt{3}$ in Regime~II arises from the
sign of the initial dual multiplier:
\[
\beta_0
=\frac{\mu t(3-L^2t^2)}{2+\mu t}\,p^{N-1}.
\]
For $Lt>\sqrt{3}$, this multiplier is negative, so that dual
certificate is no longer feasible. In contrast, the Toeplitz
primal construction remains feasible throughout $1<Lt<2$;
this will be used in
Proposition~\ref{prop:r3c2:toeplitz-full-feasible}.

\subsection{Construction of the exact optimal solution for the regime $\sqrt{3}/L <t < 2/L$}
In this section, we consider the stepsize
\[
 \sqrt{3}/L< t < 2/L.
\]
We find the solution structure has a phase transition when $N \leq  N_0(s) $ and $N> N_0(s) $ where 
\[
N_0(s) = \left\lfloor T(s) \right\rfloor= \left\lfloor \frac{\log_{s} \frac{2}{s+3}}{2} \right\rfloor.
\]
In Regime~III, optimality requires two horizon-dependent
constructions. For $N\leq N_0$, a rank-one certificate gives
the transient optimal value $c^{(N)}$. For $N>N_0$, we retain
the primal prefix of length $N_0$ and attach a geometrically
decaying Toeplitz tail. The threshold comes from the positive
semidefiniteness of the transient dual slack, not from a loss
of primal feasibility.In terms of objective gaps, in the transient phase (Theorem~\ref{thm:primal_dual_equal_case31}), primal variables $\delta_k$'s do not have a geometric decay pattern anymore. While in the geometric phase (Theorem~\ref{thm:primal_dual_equal_case32}), $\delta_k$'s return to geometric decay with rate $p$ for $k \geq N_0$.
\begin{theorem}[\bf Regime III - Phase 1]\label{thm:primal_dual_equal_case31}
Let
\[
\sqrt{3}/L < t < 2/L,
\]
and let $1 \leq N \leq N_0$ be an interger, then the following statements hold.
\begin{enumerate}
  \item \label{thm:primal_dual_equal_case31_r1}
An optimal primal solution of~\eqref{primal:reSDP} is given by
\begin{equation}\label{psol_r3c1:reSDP}
\begin{aligned}
&G_{i,j} = (-s)^{i+j-2N} 2\mu c\Delta,~0\leq j\leq i\leq N, \\
&\ell  = c\Delta, \\
&\delta_{N} = c\Delta,~\delta_{k} = \delta_{k+1} + \frac{t}{2}(G_{k,k}+G_{k,k+1}), ~0 \leq k \leq N-1
\end{aligned}
\end{equation}
where 
\begin{equation}\label{eq_r3c1:r3_cdef}
c = \frac{s^{2N}\kappa}{1+s^{2N}(\kappa-1)}.
\end{equation}
A corresponding optimal dual solution of~\eqref{Dual:reSDP} is given by
\begin{align*}
\lambda_k&= \left(\frac{1}{s+1}+\frac{1}{s^{2k+1}(s+1)} \right)c,~\lambda'_k= \left(\frac{1}{s^{2k+1}(s+1)}-\frac{s}{s+1} \right)c,~~0\leq k\leq N-1, \\
\alpha_{k} &= 0,~~0\leq k\leq N,\\
\eta_k &= 0,~~0\leq k\leq N-1,~~\eta_{N} =1,\\
\beta_k &= 0,~~0\leq k\leq N-1,~~\beta_{N}=1-c, \\
\gamma &= c.
\end{align*}
\item \label{thm:primal_dual_equal_case31_r2}
We have the solution is also an optimal solution of~\eqref{def:sdp_pep}.
\end{enumerate}
\end{theorem}
\noindent\textbf{Proof insights.}
The primal construction in the transient phase of Regime~III uses
a rank-one Gram matrix $\BG$. We impose equality in both
consecutive-pair constraints~\eqref{PP1} and~\eqref{PP2},
together with the terminal condition
\[
\delta_N=\frac{G_{N,N}}{2\mu}.
\]
These conditions lead to explicit expressions for $\BG$ and the
objective gaps $\delta_k$.

For the dual construction, we set $\alpha_k=0$ for
$0\leq k\leq N$ and $\beta_k=\eta_k=0$ for
$0\leq k\leq N-1$, leaving only the terminal multipliers
$\beta_N$ and $\eta_N$ potentially nonzero in these families.
We also choose $\lambda_k-\lambda'_k$ to be independent of $k$.
This structure simplifies the dual balance equations~\eqref{DD1}
and guides the construction of explicit dual multipliers.

Finally, we express the dual slack matrix $\BS$ as a weighted
sum of rank-one positive semidefinite matrices. Examining the
weights yields the precise condition
\[
\BS\succeq0
\quad\Longleftrightarrow\quad
N\leq N_0,
\]
see Lemma~\ref{r3:LemmaS}.
Thus, this certificate is dual feasible throughout the transient
phase and loses positive semidefiniteness beyond the threshold
$N_0$, motivating the different construction used in the
subsequent geometric phase.
\begin{proof}
The proposed primal and dual points have matching objective values, $\ell=\gamma\Delta=c\Delta$. By weak duality in Lemma~\ref{Lemma:primal_dual}, it remains only to verify all primal and dual constraints. \\
\\
\noindent\textbf{Nonnegativity.}
We first verify that all dual multipliers are nonnegative.
Since $0<s<1$, $\kappa>1$, and $N\geq1$, we have
\[
0<c
=\frac{\kappa s^{2N}}{1-s^{2N}+\kappa s^{2N}}
<1.
\]
The multipliers $\alpha_k$ and $\eta_k$ are nonnegative by
definition. Moreover, $\beta_k=0$ for $0\leq k\leq N-1$,
$\beta_N=1-c>0$, and $\gamma=c>0$.
For $0\leq k\leq N-1$, the expression for $\lambda_k$
is positive, while
\[
\lambda'_k
=\frac{c}{s+1}\left(s^{-2k-1}-s\right)
\geq\frac{c(1-s)}{s+1}
>0.
\]
Thus all dual multipliers are nonnegative.\\
\\
\textbf{Primal Feasibility.} We verify primal constraints~\eqref{PP1}--\eqref{PP6} in order.\\
\\
\eqref{PP1},\eqref{PP2}: These constraints hold with equality by Lemma~\ref{lem:r3c1:full-primal}.\\
\\
\eqref{PP3}: For $k=N$, we have $\delta_{N} = G_{N,N}/(2\mu)$. We start by calculating the expression for $\delta_{k},~~0\leq k\leq N-1$, we have
\begin{equation}\label{eq:R3C1:EQ3}
\begin{aligned}
\delta_{k} &= \delta_{N} + \sum_{i=k}^{N-1} (\delta_{i} - \delta_{i+1}) \\
&=\delta_{N} + \frac{t}{2} \sum_{i=k}^{N-1} \left(G_{i,i}+G_{i,i+1}\right) \\
&=c\Delta + \frac{t(1-s) }{2}\sum_{i=k}^{N-1} G_{i,i} \\
&=c\Delta  + \mu t (1-s) c\Delta \sum_{i=k}^{N-1} s^{2i-2N} \\
& = c\Delta + \mu t(1-s)c\Delta \frac{1 - s^{2k-2N}}{s^2-1} \\
& = c\Delta \left(1+ \kappa^{-1}(s^{2k-2N}-1)\right)
\end{aligned}
\end{equation}
where
\[ 
\mu t(1-s) = \kappa^{-1}(1+s)(1-s) = \kappa^{-1}(1-s^2),~~~~
\sum_{i=k}^{N-1} s^{2i-2N}  = \frac{1 - s^{2k-2N} }{s^2-1}.
\]
Hence we have
\[
\frac{1}{2L} \leq \frac{\delta_{k}}{G_{k,k}} = \left(\left(1-\frac{1}{\kappa} \right)s^{2N-2k} + \frac{1}{\kappa}\right)\frac{1}{2\mu} = \frac{1}{2L} + \left(\frac{1}{2\mu}-\frac{1}{2L} \right) s^{2N-2k} \leq \frac{1}{2\mu}.
\]
as $0< s < 1$ and $k \leq N-1$.
Now we have verified the constraints. \\
\\
\eqref{PP4}: By~\eqref{eq:R3C1:EQ3}, we have $\delta_0= c\Delta \left(1+ \kappa^{-1}(s^{-2N}-1)\right)=c\Delta c^{-1} = \Delta$.\\
\\
\eqref{PP5}: By~\eqref{eq:R3C1:EQ3}, we have $\ell = \delta_{N} <\delta_{N-1} < \cdots<\delta_0$. \\
\\
\eqref{PP6}: We have by definition
\[
\BG = 2\mu c\Delta \bv\bv^{\top} \succeq 0,~~v_k = (-s)^{k-N}.
\]
\textbf{Dual Feasibility.}  We proceed to verify the dual constraints ~\eqref{DD1}--\eqref{DD3} in order. \\
\\
\noindent\eqref{DD1}:
The proposed multipliers satisfy
\[
\lambda_k-\lambda'_k=c,
\qquad 0\leq k\leq N-1.
\]
At the initial index $k=0$, we have
$\alpha_0=\beta_0=\eta_0=0$ and $\gamma=c$, so
\[
\beta_0-\alpha_0-\eta_0-\lambda_0+\lambda'_0+\gamma
=-(\lambda_0-\lambda'_0)+c
=-c+c=0.
\]
For each interior index $1\leq k\leq N-1$, the multipliers
$\alpha_k$, $\beta_k$, and $\eta_k$ vanish. Hence
\begin{align*}
&\beta_k-\alpha_k-\eta_k
+\lambda_{k-1}-\lambda_k-\lambda'_{k-1}+\lambda'_k\\
&\qquad
=(\lambda_{k-1}-\lambda'_{k-1})
-(\lambda_k-\lambda'_k)
=c-c=0.
\end{align*}
At the terminal index $k=N$, we have
$\alpha_N=0$, $\beta_N=1-c$, and $\eta_N=1$, giving
\[
\beta_N-\alpha_N-\eta_N+\lambda_{N-1}-\lambda'_{N-1}
=(1-c)-1+c=0.
\]
Thus all balance equations in~\eqref{DD1} hold.\\
\\
\eqref{DD2}: The constraint holds trivially since $\eta_{N}= 1$ and $\eta_{k}=0$ for all $0\leq k\leq N-1$.\\
\\
\eqref{DD3}: The constraint holds by Lemma~\ref{r3:LemmaS} since $N\leq N_0=\lfloor T(s)\rfloor\leq T(s)$. \\
\\
\textbf{Full-SDP Feasibility.} We proceed with item 2 in Theorem~\ref{thm:primal_dual_equal_case31}. Since~\eqref{primal:reSDP} is a relaxation of~\eqref{def:sdp_pep}, It remains to verify $\BG$ and $\left\{\delta_k\right\}$ are feasible for non-consecutive pairs for the first constraint in~\eqref{def:sdp_pep}. This condition is directly verified Lemma~\ref{lem:r3c1:full-primal}. Extend the reduced dual certificate by assigning zero multipliers
to the omitted nonconsecutive constraints. The resulting full-SDP
primal and dual points remain feasible and have matching objective
values. Weak duality therefore establishes their optimality for
the full pairwise SDP.
\end{proof}

\begin{theorem}[\bf Regime III - Phase 2]\label{thm:primal_dual_equal_case32}
Let
\[
\sqrt{3}/L < t < 2/L,
\]
and let $N>N_0$ be an integer, then the following statements hold.
\begin{enumerate}
\item \label{thm:primal_dual_equal_case3_p1_r1} An optimal primal solution to~\eqref{primal:reSDP} is given by
\begin{equation}\label{psol_r3c2:reSDP}
\begin{aligned}
\ell &= p^{N-N_0}c^{(N_0)}\Delta, \\
\delta_k & = 
\begin{dcases}
c^{(N_0)} \Delta, & k = N_0, \\
\delta_{k+1} + \frac{t}{2}(G_{k,k}+G_{k,k+1}), & 0 \leq k\leq N_0-1, \\
p^{k-N_0}\delta_{N_0},  &  N_0 <  k \leq N,
\end{dcases} \\
%
G_{i,j}  & = 2\mu c^{(N_0)}\Delta
\begin{dcases}
(-s)^{i+j-2N_0}, & 0\leq i, j \leq  N_0-1, \\
(-s)^{i-N_0} p^{(j-N_0)/2}\phi_{j-N_0}, & 0 \leq i < N_0 \leq j \leq N,~~~~G_{j,i} = G_{i,j}  \\
 [\BP^{(N-N_0)} \BPhi^{(N-N_0)}\BP^{(N-N_0)}]_{i-N_0,j-N_0}, & N_0 \leq i,j\leq N,
\end{dcases}
\end{aligned}
\end{equation}
where 
\begin{equation}\label{eq_r3c2:r3_cdef}
c^{(N)} = \frac{s^{2N}\kappa}{1+s^{2N}(\kappa-1)} = \frac{(Lt-1)^{2N}L}{\mu+(Lt-1)^{2N}(L-\mu)}.
\end{equation}
and $\BPhi$ and $\BP$ are defined in \eqref{def:BPBPHI} and \eqref{R2:eq_recur}. The corresponding dual variables to~\eqref{Dual:reSDP}
are given by
\begin{align*}
\gamma
&=p^{N-N_0}c^{(N_0)},\\
\alpha_k
&=0,\qquad 0\leq k\leq N,\\
\eta_k
&=0,\qquad 0\leq k<N,
\qquad \eta_N=1,\\
\beta_k
&=
\begin{cases}
0,
& 0\leq k<N_0,\\[1mm]
\displaystyle
p^{N-N_0-1}
\left(
\frac{2}{2+\mu t}-p\,c^{(N_0)}
\right),
& k=N_0,\\[3mm]
\displaystyle
\frac{2(1-p)p^{N-k-1}}{2+\mu t},
& N_0<k<N,\\[3mm]
\displaystyle
\frac{\mu t}{2+\mu t},
& k=N,
\end{cases}\\
\lambda_k
&=
\begin{cases}
p^{N-N_0}\lambda_k^{(N_0)},
& 0\leq k<N_0,\\[1mm]
\displaystyle
\frac{(s+2)p^{N-k-1}}{2+\mu t},
& N_0\leq k<N,
\end{cases}\\
\lambda'_k
&=
\begin{cases}
p^{N-N_0}\lambda_k^{\prime(N_0)},
& 0\leq k<N_0,\\[1mm]
\displaystyle
\frac{s\,p^{N-k-1}}{2+\mu t},
& N_0\leq k<N.
\end{cases}
\end{align*}
Here $\lambda_k^{(N_0)}$ and
$\lambda_k^{\prime(N_0)}$ are the dual multipliers from
Theorem~\ref{thm:primal_dual_equal_case31}, evaluated at
horizon $N_0$. When $N=N_0+1$, the branch
$N_0<k<N$ in the definition of $\beta_k$ is empty.
\item \label{thm:primal_dual_equal_case3_p1_r2} We have the solution is not a feasible point of~\eqref{def:sdp_pep}.
\end{enumerate}
\end{theorem}
\noindent\textbf{Proof insights.}
The key observation is that the optimal solution combines the
transient construction from Phase~1 with the geometric structure
of Regime~II. The primal variables on indices $0,\ldots,N_0$
coincide with the horizon-$N_0$ transient solution. Beyond this
prefix, the objective gaps and the diagonal entries of $\BG$
decay geometrically with ratio $p$.

For the primal construction, we choose
$\BG_{N_0:N,N_0:N}$ to be the Toeplitz-envelope matrix used in
Regime~II, evaluated at the current stepsize and scaled by
$c^{(N_0)}$. The prefix and tail overlap at index $N_0$, where
\[
G_{N_0,N_0}=2\mu c^{(N_0)}\Delta
           =2\mu\delta_{N_0}.
\]
To complete the Gram matrix, we choose the cross-block entries
so that each prefix gradient vector remains a scalar multiple
of the gradient vector at the transition. This yields the
factorization in~\eqref{eq:form_G}, which establishes
$\BG\succeq0$ from the positive semidefiniteness of the
Toeplitz block. Unlike the Regime-II construction, the full
mixed Gram matrix is singular.

For the dual construction, we scale the Phase-1 multipliers
$\lambda_k^{(N_0)}$ and $\lambda_k^{\prime(N_0)}$ by
$p^{N-N_0}$ for $k<N_0$, and use the geometric formulas from
Regime~II for $k\geq N_0$. The dual balance equation~\eqref{DD1}
at the transition determines
\[
\beta_{N_0}
=
p^{N-N_0-1}
\left(\frac{2}{2+\mu t}-p\,c^{(N_0)}\right).
\]
Its positivity follows from the threshold relation
$N_0+1>T(s)$. The interior tail multipliers then follow the
same geometric pattern as in Regime~II, while the terminal
multipliers satisfy the final balance equation.

Finally, the tail contributions to the dual slack vanish,
and the transition multiplier makes the boundary diagonal
entry agree with the scaled transient certificate. Consequently,
\[
\BS^{(N)}
=
p^{N-N_0}
\begin{bmatrix}
\BS^{(N_0)} & 0\\
0 & 0
\end{bmatrix}
\succeq0.
\]
Thus dual feasibility follows from the horizon-$N_0$
Phase-1 certificate. Together with primal feasibility and
the matching objective values
$\ell=\gamma\Delta=c^{(N_0)}p^{N-N_0}\Delta$,
this establishes optimality for the consecutive-pair SDP.
\begin{proof}
The proposed primal and dual points have matching objective values, $\ell = \gamma\Delta = p^{N-N_0}c^{(N_0)}\Delta$. By weak duality in Lemma~\ref{Lemma:primal_dual}, it remains only to verify all primal and dual constraints. \\
\\
\noindent\textbf{Nonnegativity.}
We first verify that all dual multipliers are nonnegative.
By Lemma~\ref{R2:L1}, we have $0<p<1$.
Moreover, $0<s<1$ and $c^{(N_0)}>0$.
The multipliers $\alpha_k$ and $\eta_k$ are nonnegative by
definition, while
\[
\gamma=p^{N-N_0}c^{(N_0)}>0.
\]
\\
For $0\leq k<N_0$, the multipliers $\lambda_k$ and $\lambda'_k$
are positive rescalings of the Phase-1 multipliers:
\[
\lambda_k=p^{N-N_0}\lambda_k^{(N_0)}>0,
\qquad
\lambda'_k=p^{N-N_0}\lambda_k^{\prime(N_0)}>0,
\]
where positivity follows from
Theorem~\ref{thm:primal_dual_equal_case31}.
For $N_0\leq k\leq N-1$, we have
$$
\lambda_k=\frac{s+2}{2+\mu t}\,p^{N-k-1}>0,~~\lambda'_k=\frac{s}{2+\mu t}\,p^{N-k-1}>0.
$$
\\
For the multipliers $\beta_k$ away from the transition index,
we obtain
\begin{align*}
\beta_k&=0,~~
0\leq k<N_0,\\
\beta_k
&=\frac{2(1-p)}{2+\mu t}\,p^{N-k-1}>0,~~
N_0<k<N,\\
\beta_N&=\frac{\mu t}{2+\mu t}>0.
\end{align*}
When $N=N_0+1$, the interior range $N_0<k<N$ is empty. It remains to verify the sign of the boundary multiplier
$\beta_{N_0}$.
Using
\[
p=\frac{s^3+3s^2+2\kappa-2}{2\kappa+s+1},
\qquad
c^{(N_0)}
=\frac{\kappa s^{2N_0}}{1+(\kappa-1)s^{2N_0}},
\qquad
\mu t=\frac{s+1}{\kappa},
\]
we obtain
\begin{align*}
\beta_{N_0}
&=p^{N-N_0-1}
  \left(\frac{2}{2+\mu t}-p\,c^{(N_0)}\right) =\frac{\kappa p^{N-N_0-1}\left[2-(s+3)s^{2N_0+2}\right]}{(2\kappa+s+1)\left[1+(\kappa-1)s^{2N_0}\right] }.
\end{align*}
By the definition of the transition threshold,
\[
N_0=\lfloor T(s)\rfloor, \qquad T(s):=\frac{\log\!\left(2/(s+3)\right)}{2\log s},
\]
so $N_0+1>T(s)$.
Since $0<s<1$, it follows that
\[
s^{2N_0+2}<s^{2T(s)}=\frac{2}{s+3},
\]
and hence
\[
2-(s+3)s^{2N_0+2}>0.
\]
All remaining factors in the expression for $\beta_{N_0}$
are positive, so $\beta_{N_0}>0$.
Thus all dual multipliers are nonnegative.\\
\\
\textbf{Primal Feasibility:} We verify all the primal constraints \eqref{PP1}--\eqref{PP6} in order. \\
\\
\eqref{PP1},\eqref{PP2}: For any $0\leq k\leq N_0-1$, we have $\delta_{k}^{(N)}  = \delta_k^{(N_0)} $ and $\BG_{0:N_0,0:N_0}^{(N)} = \BG^{(N_0)}$, hence the constraints are verified by Theorem~\ref{thm:primal_dual_equal_case31}. For any $ N_0 \leq k \leq N-1$,~we have $\delta_{k}$ and $\BG_{N_0:N,N_0:N}$ have the same structure as primal variables in Theorem~\ref{thm:primal_dual_equal_case2} when number of steps is given by $N-N_0$ scaled by $c^{(N_0)}$ these constraints holds by Lemma~\ref{R2:L2} as this lemma holds for $1/L<t<2/L$.\\
\\
~\eqref{PP3}:  For $k \leq N_0$, since variables are identical to $N=N_0$, it has been verified in~Theorem~\ref{thm:primal_dual_equal_case31}. For $k>N_0$, we have by definition
\[
\delta_k = \frac{1}{2\mu} G_{k,k} \geq \frac{1}{2L}G_{k,k}.
\]
\\
\eqref{PP4}: By definition we have $\delta_{0} = \Delta$. \\
\\
\eqref{PP5}: We have $p<1$ by Lemma~\ref{R2:L1}, hence
\[
\ell = \delta_{N} <\delta_{N-1} < \cdots<\delta_{N_0}.
\]
Additionally, by~\eqref{eq:R3C1:EQ3} it holds that
\[
\delta_{N_0} < \delta_{N_0-1} < \cdots < \delta_0.
\]
These two equations allow us to conclude $\delta_{k} \geq \ell$ for all $k$. \\
\\
\eqref{PP6}: We have 
\begin{align}
\BG &=2\mu c^{(N_0)}\Delta \begin{bmatrix}\bv \bv^{\top}  & \bv (\be^{(N-N_0)}_0)^{\top}\BB^{(N-N_0)}\\
\BB^{(N-N_0)} \be^{(N-N_0)}_0 \bv^{\top}  & \BB^{(N-N_0)} \end{bmatrix},\\
v_{k}:&= (-s)^{k-N_0}~(k < N_0),~\BB^{(N-N_0)}:= \BP^{(N-N_0)}\BPhi^{(N-N_0)}\BP^{(N-N_0)}.
\end{align}
Hence we have
\begin{equation}\label{eq:form_G}
\BG = 2\mu c^{(N_0)}\Delta \begin{bmatrix}\bv (\be^{(N-N_0)}_0)^{\top} \\ \I_{N-N_0+1} \end{bmatrix}\BB^{(N-N_0)}\begin{bmatrix}\be^{(N-N_0)}_0 \bv^{\top} & \I_{N-N_0+1} \end{bmatrix} \succeq 0,
\end{equation}
Since by Lemma~\ref{R2:L3}, $\BB^{(N-N_0)} \succeq 0$. \\
\\
\textbf{Dual Feasibility.} We proceed to verify the dual constraints \eqref{DD1}--\eqref{DD3} in order.\\
\\
\eqref{DD1}: The balance equations follow for $k < N_0$ since the variables are identical to the case when $N=N_0$ scaled by $p^{N-N_0}$. For $k = N_0$, we have
\[
(\lambda_{N_0} - \lambda'_{N_0}) + (\lambda'_{N_0-1}  - \lambda_{N_0-1}) = p^{N-N_0-1}\frac{2}{2+\mu t} - p^{N-N_0}c^{(N_0)} = p^{N-N_0-1}\left(\frac{2}{2+\mu t} - pc^{(N_0)}\right) = \beta_{N_0}.
\]
For $N_0 < k < N$, we have
\[
(\lambda_{k} - \lambda'_{k}) + (\lambda'_{k-1}  - \lambda_{k-1}) =p^{N-k-1} \frac{2}{2+\mu t} - p^{N-k} \frac{2}{2+\mu t} = p^{N-k-1}\frac{2(1-p)}{2+\mu t} = \beta_k.
\]
For $k = N$, we have
\[
\beta_{N}-\eta_{N} = \frac{\mu t}{2+\mu t}-1 = -\frac{2}{2+\mu t} = \lambda_{N-1}' - \lambda_{N-1}.
\]
\\      
\eqref{DD2}: The constraint holds trivially since $\eta_{N}= 1$ and $\eta_{k}=0,~0 \leq k \leq N-1$.\\
\\
\eqref{DD3}: Notice for block $\BS_{N_0+1:N,N_0+1:N}$, the dual variables structure is identical to the regime-$2$ case with
$$
\lambda_{k} = \frac{s+2}{s} \lambda_{k}',~~\beta_k = \lambda_{k-1}'-\lambda_{k}'+\lambda_{k} - \lambda_{k-1},~~N_0 \leq k \leq N-1,~~\beta_{N}-1=\lambda'_{N-1}-\lambda_{N-1}.
$$
Hence from the same calculation in Lemma~\ref{R2:L4}, we have $\BS_{N_0+1:N,N_0+1:N} = 0$. For block $\BS_{0:N_0-1,0:N_0-1}$, we have 
$$
\BS^{(N)}_{0:N_0-1,0:N_0-1} = p^{N-N_0}\BS^{(N_0)}_{0:N_0-1,0:N_0-1},~~S^{(N)}_{N_0-1,N_0} = S^{(N)}_{N_0,N_0-1}=p^{N-N_0}S^{(N_0)}_{N_0-1,N_0}.
$$
Now we compute
\begin{align*}
 A_{0,0}\lambda_{N_0}+ A'_{0,0}\lambda'_{N_0} &= \frac{p^{N-N_0-1}}{4L(2+\mu t)} \left[(2-s^2)(s+2)-(s^2+4s+2)s\right] \\
 & = \frac{p^{N-N_0-1}}{2L(2+\mu t)} \left(-s^3-3s^2 + 2\right) \\
 & = \frac{p^{N-N_0-1}}{2(2+\mu t)} \left(3t-L^2t^3 \right) \\
 & = \frac{p^{N-N_0-1}}{2\mu}\left(\frac{2}{2+\mu t}+\frac{3\mu t-\mu L^2t^3-2}{2+\mu t} \right) \\
 & = \frac{p^{N-N_0-1}}{2\mu}\left(\frac{2}{2+\mu t}-p\right)
\end{align*}
and
$$
\beta_{N_0}^{(N)} - p^{N-N_0}\beta_{N_0}^{(N_0)} = p^{N-N_0-1}\left(\frac{2}{2+\mu t}-p\,c^{(N_0)}\right) -  p^{N-N_0-1}(p-pc^{(N_0)}) =  p^{N-N_0-1}\left(\frac{2}{2+\mu t}-p\right).
$$
Hence we have $2\mu(A_{0,0}\lambda_{N_0}+ A'_{0,0}\lambda'_{N_0}) = \beta_{N_0}^{(N)} - p^{N-N_0}\beta_{N_0}^{(N_0)} $ and we can compute
\begin{align*}
S^{(N)}_{N_0,N_0} &=A_{1,1}\lambda^{(N)}_{N_0-1}  + A'_{1,1}(\lambda')^{(N)}_{N_0-1}
  + A_{0,0}\lambda^{(N)}_{N_0}+ A'_{0,0}(\lambda')^{(N)}_{N_0} -\frac{1}{2\mu}\beta^{(N)}_{N_0} \\
  &= A_{1,1}\lambda^{(N)}_{N_0-1}  + A'_{1,1}(\lambda')^{(N)}_{N_0-1} + \frac{1}{2\mu}\left(\beta_{N_0}^{(N)} - p^{N-N_0}\beta_{N_0}^{(N_0)}\right) -\frac{1}{2\mu}\beta^{(N)}_{N_0} \\
  &= p^{N-N_0}\left(A_{1,1}\lambda^{(N_0)}_{N_0-1}  + A'_{1,1}(\lambda')^{(N_0)}_{N_0-1} -\frac{1}{2\mu}\beta^{(N_0)}_{N_0}\right) \\
  &= p^{N-N_0}S^{(N_0)}_{N_0,N_0}.
\end{align*}
The cross-boundary off-diagonal entries also vanish:
\[
S^{(N)}_{N_0,N_0+1}
=S^{(N)}_{N_0+1,N_0}
=\frac{s\lambda_{N_0}-(s+2)\lambda'_{N_0}}{4L}
=0.
\]
Since $\BS^{(N)}$ is tridiagonal, the preceding calculations
establish the claimed block-diagonal identity.
This allows us to conclude that 
$$
\BS^{(N)} = p^{N-N_0}\begin{bmatrix} \BS^{(N_0)} & \bzero_{N_0+1,N-N_0} \\  \bzero_{N-N_0,N_0+1} &  \bzero_{N-N_0,N-N_0} \end{bmatrix} \succeq 0
$$
by \eqref{DD3} in Theorem~\ref{thm:primal_dual_equal_case31}. Hence we have verified all dual constraints. \\
\\
\textbf{Infeasibility of the full-SDP.} The primal solution is not feasible for the full-SDP~\eqref{def:sdp_pep} because it does not satisfy nonconsecutive smoothness interpolation inequalities. In particular, for the pair $(N_0+1,N_0-1)$, we have
$$
\Tr(\BA_{N_0+1,N_0-1}\BG) - (\delta_{N_0+1}-\delta_{N_0-1}) = \frac{2c\Delta (\kappa-1)(1-s)(1+s)^2}{\kappa s(2\kappa+s+1)} > 0.
$$
\end{proof}
\medskip
\noindent\textbf{A full-SDP feasible construction.}
Theorem~\ref{thm:primal_dual_equal_case32} determines the exact
value of the consecutive-pair SDP for $N>N_0$. However,
item~\ref{thm:primal_dual_equal_case3_p1_r2} shows that the
constructed mixed primal optimizer is not feasible for the
full pairwise SDP~\eqref{def:sdp_pep}. Thus, this construction
does not establish finite-horizon equality between the two
SDP values, nor does its infeasibility imply a strict gap
between them.

To identify the asymptotic factor of the full pairwise SDP,
we instead extend the diagonally scaled Toeplitz construction
from Regime~II over the entire horizon, without a rank-one
prefix. The following proposition establishes that this
alternative point is full-SDP feasible and has normalized
objective value $p(t)^N$. Although it need not attain the
finite-horizon optimum, it provides a lower bound that,
together with the reduced-SDP upper bound, determines the
exact asymptotic geometric factor of the full pairwise SDP.

\begin{proposition}[A full-SDP feasible point in the geometric phase]
\label{prop:r3c2:toeplitz-full-feasible}
For
\[
\sqrt{3}/L < t < 2/L,
\]
and let $N>N_0$ be an integer and
\[
 p:=\frac{\mu L^2t^3-3\mu t+2}{2+\mu t},
 \qquad \tau:=\frac{1-L^2t^2+p}{2}.
\]
We define
\begin{equation}\label{eq:r3c2:toeplitz-feasible-point}
 \delta_k=\Delta p^k,~~0\leq k\leq N,
 \qquad \ell=\Delta p^N,
 \qquad \BG=2\mu\Delta\BP^{(N)}\BPhi^{(N)}\BP^{(N)},
\end{equation}
where $\BP$ and $\BPhi$ are from~\eqref{def:BPBPHI} with the recurrence
\eqref{R2_proof:eq_recur} evaluated at these values of $s$ and $p$.
Then~\eqref{eq:r3c2:toeplitz-feasible-point} is feasible for the full
pairwise SDP~\eqref{def:sdp_pep}. Consequently,
\begin{equation}\label{eq:r3c2:full-sdp-lower-bound}
 p(t)^N\leq\overline{\mathcal R}_N(t).
\end{equation}
This is the Regime-II primal construction extended over the entire
horizon $0,\ldots,N$, not the mixed optimizer obtained by joining a
rank-one prefix to a Toeplitz tail.
\end{proposition}

\begin{proof}
The primal-feasibility lemmas in the Regime-II appendix are stated
for the larger interval $1/L < t<2/L$. Lemma~\ref{R2:L3} gives $\BG\succ0$, and
Lemma~\ref{R2:L2} verifies every ordered-pair interpolation
constraint, including all nonconsecutive pairs. Additionally, we directly plug in to verify the remaining primal constraints:
$$
\delta_0 = \Delta,~~\ell=\delta_{N} \leq \delta_k,~~G_{k,k} = 2\mu\delta_k,~0 \leq k \leq N.
$$
Thus~\eqref{eq:r3c2:toeplitz-feasible-point} is full-SDP primal
feasible. Its objective value is $\Delta p^N$, whereas the
full-SDP optimal value is $\Delta\overline{\mathcal R}_N(t)$.
Maximization and division by $\Delta>0$ yield
\eqref{eq:r3c2:full-sdp-lower-bound}.
\end{proof}

Combining
Proposition~\ref{prop:r3c2:toeplitz-full-feasible}
with Theorem~\ref{thm:primal_dual_equal_case32},
for every $N>N_0$ we obtain
\[
p(t)^N
\leq\overline{\mathcal R}_N(t)
\leq\widehat{\mathcal R}_N(t)
=c^{(N_0)}p(t)^{N-N_0}.
\]
For each fixed $t\in(\sqrt{3}/L,2/L)$, the quantities
$c^{(N_0)}$ and $N_0$ are independent of $N$.
Taking $N$th roots and passing to the limit gives
\[
\lim_{N\to\infty}\overline{\mathcal R}_N(t)^{1/N}
=
\lim_{N\to\infty}\widehat{\mathcal R}_N(t)^{1/N}
=p(t).
\]
This identifies the exact asymptotic factor of the full
pairwise SDP without asserting finite-horizon equality in
the geometric phase. The feasible SDP point does not supply
a lower bound for the functional PEP unless it can also be
realized by a globally smooth PL function.

\section{Concluding remarks}\label{s:conclusion}

We have derived explicit multi-step convergence bounds for
constant-stepsize gradient descent on $L$-smooth functions satisfying
the $\mu$-P\L{} inequality. Matching primal--dual certificates solve
the consecutive-pair SDP analytically for every horizon $N\geq1$
and every stepsize $0<t<2/L$. The analysis identifies exact
tensorization of the one-step SDP factor for $0<t\leq\sqrt{3}/L$
and a nongeometric transient followed by a geometric tail for
$\sqrt{3}/L<t<2/L$. In the latter regime, the multi-step analysis
strictly improves the asymptotic factor obtained by recursively
applying the one-step estimate of~\cite{AKZ23}.

These bounds are the sharpest obtainable within the consecutive-pair
SDP. They also attain the full pairwise SDP value in Regimes~I and~II
and during the transient phase of Regime~III, while their asymptotic
geometric factor is exact for the full pairwise SDP in all regimes.
Consequently, no smaller asymptotic factor can be derived from the
same pairwise smoothness and pointwise P\L{} constraints alone.
The full-SDP finite-horizon optimum in the geometric phase of
Regime~III remains to be determined, and SDP exactness does not
establish tightness over the underlying function class.

A central direction for further progress is an exact interpolation
characterization for smooth P\L{} functions: necessary and sufficient
conditions on finite samples of function values and gradients that
ensure the existence of a global interpolant with the prescribed
constants $L$ and $\mu$. Such a characterization would clarify the
gap between the present SDP relaxations and the functional PEP,
provide a systematic route to matching worst-case functions, and
guide further improvements beyond the current framework.

\bibliographystyle{abbrv}
\bibliography{ref}

@article{DZPS19,
	author = {Du, Simon S and Zhai, Xiyu and Poczos, Barnabas and Singh, Aarti},
	journal = {ICLR},
	title = {Gradient descent provably optimizes over-parameterized neural networks},
	year = {2019}}

@article{TV23,
	author = {Teboulle, Marc and Vaisbourd, Yakov},
	doi = {10.1007/s10107-022-01899-0},
	journal = {Mathematical Programming},
	number = {1--2},
	pages = {63--96},
	title = {An Elementary Approach to Tight Worst-Case Complexity Analysis of Gradient-Based Methods},
	volume = {201},
	year = {2023}}

@article{AP25a,
	articleno = {12},
	author = {Altschuler, Jason M. and Parrilo, Pablo A.},
	doi = {10.1145/3708502},
	journal = {Journal of the ACM},
	number = {2},
	numpages = {38},
	pages = {1--38},
	title = {Acceleration by Stepsize Hedging: Multi-Step Descent and the Silver Stepsize Schedule},
	volume = {72},
	year = {2025}}

@article{AP25b,
	author = {Altschuler, Jason M. and Parrilo, Pablo A.},
	journal = {Mathematical Programming},
	number = {1},
	pages = {1105--1118},
	title = {Acceleration by Stepsize Hedging: Silver Stepsize Schedule for Smooth Convex Optimization},
	volume = {213},
	year = {2025}}

@article{G24,
	author = {Grimmer, Benjamin},
	journal = {SIAM Journal on Optimization},
	number = {3},
	pages = {2588--2608},
	title = {Provably Faster Gradient Descent via Long Steps},
	volume = {34},
	year = {2024}}

@article{GSW25,
	author = {Grimmer, Benjamin and Shu, Kevin and Wang, Alex L.},
	doi = {10.1287/ijoo.2024.0057},
	journal = {INFORMS Journal on Optimization},
	number = {2},
	pages = {156--169},
	title = {Accelerated Objective Gap and Gradient Norm Convergence for Gradient Descent via Long Steps},
	volume = {7},
	year = {2025}}

@article{BI95,
	author = {Bonnans, Joseph Fr{\'e}d{\'e}ric and Ioffe, Alexander},
	journal = {Mathematics of Operations Research},
	number = {4},
	pages = {801--817},
	publisher = {INFORMS},
	title = {Second-order sufficiency and quadratic growth for nonisolated minima},
	volume = {20},
	year = {1995}}

@article{M00,
	author = {Anitescu, Mihai},
	journal = {SIAM Journal on Optimization},
	number = {4},
	pages = {1116--1135},
	publisher = {SIAM},
	title = {Degenerate nonlinear programming with a quadratic growth condition},
	volume = {10},
	year = {2000}}

@article{LT93,
	author = {Luo, Zhi-Quan and Tseng, Paul},
	journal = {Annals of Operations Research},
	number = {1},
	pages = {157--178},
	publisher = {Springer},
	title = {Error bounds and convergence analysis of feasible descent methods: a general approach},
	volume = {46},
	year = {1993}}

@inproceedings{Loj63,
	address = {Paris},
	author = {{\L}ojasiewicz, Stanis{\l}aw},
	booktitle = {Les {\'E}quations aux D{\'e}riv{\'e}es Partielles (Paris, 1962)},
	mrnumber = {MR0160856},
	number = {117},
	pages = {87--89},
	publisher = {{\'E}ditions du Centre National de la Recherche Scientifique},
	series = {Colloques internationaux du Centre national de la recherche scientifique},
	title = {Une propri{\'e}t{\'e} topologique des sous-ensembles analytiques r{\'e}els},
	year = {1963}}

@article{AKZ22,
	author = {Abbaszadehpeivasti, Hadi and de Klerk, Etienne and Zamani, Moslem},
	journal = {Optimization Letters},
	number = {6},
	pages = {1649--1661},
	publisher = {Springer},
	title = {The exact worst-case convergence rate of the gradient method with fixed step lengths for $L$-smooth functions},
	volume = {16},
	year = {2022}}

@article{KJH26,
	author = {Kong, Qingyuan and Jiang, Rujun and He, Yihan},
	journal = {Optimization Letters},
	number = {3},
	pages = {537--557},
	publisher = {Springer},
	title = {Linear convergence of the proximal gradient method for composite optimization under the {Polyak-\L ojasiewicz} inequality and its variant},
	volume = {20},
	year = {2026}}

@inproceedings{YFL23,
	author = {Yue, Pengyun and Fang, Cong and Lin, Zhouchen},
	booktitle = {The Thirty Sixth Annual Conference on Learning Theory},
	organization = {PMLR},
	pages = {2948--2968},
	title = {On the lower bound of minimizing polyak-\L ojasiewicz functions},
	year = {2023}}

@article{KK26,
	author = {Kassing, Sebastian and Kruse, Thomas},
	journal = {arXiv preprint arXiv:2605.14663},
	title = {Optimal Asymptotic Rates for (Stochastic) Gradient Descent under the Local PL-Condition: A Geometric Approach},
	year = {2026}}

@article{RGP26,
	author = {Rotaru, Teodor and Glineur, Fran{\c{c}}ois and Patrinos, Panagiotis},
	journal = {Mathematical Programming},
	pages = {1--67},
	publisher = {Springer},
	title = {Exact worst-case convergence rates of gradient descent: a complete analysis for all constant stepsizes over nonconvex and convex functions},
	year = {2026}}

@article{SL16,
	author = {Sun, Ruoyu and Luo, Zhi-Quan},
	journal = {IEEE Transactions on Information Theory},
	number = {11},
	pages = {6535--6579},
	publisher = {IEEE},
	title = {Guaranteed matrix completion via non-convex factorization},
	volume = {62},
	year = {2016}}

@article{LLSW19,
	author = {Li, Xiaodong and Ling, Shuyang and Strohmer, Thomas and Wei, Ke},
	journal = {Applied and Computational Harmonic Analysis},
	number = {3},
	pages = {893--934},
	publisher = {Elsevier},
	title = {Rapid, robust, and reliable blind deconvolution via nonconvex optimization},
	volume = {47},
	year = {2019}}

@article{LZB22,
	author = {Liu, Chaoyue and Zhu, Libin and Belkin, Mikhail},
	journal = {Applied and Computational Harmonic Analysis},
	pages = {85--116},
	publisher = {Elsevier},
	title = {Loss landscapes and optimization in over-parameterized non-linear systems and neural networks},
	volume = {59},
	year = {2022}}

@article{RHT26,
	author = {Rubbens, Anne and Hendrickx, Julien M and Taylor, Adrien},
	journal = {Journal of Nonsmooth Analysis and Optimization},
	publisher = {Episciences. org},
	title = {A constructive approach to strengthen algebraic descriptions of function and operator classes},
	volume = {6},
	year = {2026}}

@article{DVR24,
	author = {Das Gupta, Shuvomoy and Van Parys, Bart PG and Ryu, Ernest K},
	journal = {Mathematical Programming},
	number = {1},
	pages = {567--639},
	publisher = {Springer},
	title = {Branch-and-bound performance estimation programming: a unified methodology for constructing optimal optimization methods},
	volume = {204},
	year = {2024}}

@article{BCR26,
	author = {Boumal, Nicolas and Criscitiello, Christopher and Rebjock, Quentin},
	journal = {arXiv preprint arXiv:2604.07972},
	title = {Smooth, globally {Polyak-{\L}ojasiewicz} functions are nonlinear least-squares},
	year = {2026}}

@article{RB25,
	author = {Rebjock, Quentin and Boumal, Nicolas},
	journal = {Mathematical Programming},
	number = {1},
	pages = {151--199},
	publisher = {Springer},
	title = {Fast convergence to non-isolated minima: four equivalent conditions for $\mathcal{C}^2$ functions},
	volume = {213},
	year = {2025}}

@article{BCF22,
	author = {Goujaud, Baptiste and Moucer, C{\'e}line and Glineur, Fran{\c{c}}ois and Hendrickx, Julien M and Taylor, Adrien B and Dieuleveut, Aymeric},
	journal = {Mathematical Programming Computation},
	number = {3},
	pages = {337--367},
	publisher = {Springer},
	title = {{PEPit}: computer-assisted worst-case analyses of first-order optimization methods in Python},
	volume = {16},
	year = {2024}}

@article{AKZ23,
	author = {Abbaszadehpeivasti, Hadi and de Klerk, Etienne and Zamani, Moslem},
	journal = {Optimization Letters},
	number = {5},
	pages = {1105--1125},
	publisher = {Springer},
	title = {Conditions for linear convergence of the gradient method for non-convex optimization},
	volume = {17},
	year = {2023}}

@article{THG17,
	author = {Taylor, Adrien B and Hendrickx, Julien M and Glineur, Fran{\c{c}}ois},
	journal = {Mathematical Programming},
	number = {1},
	pages = {307--345},
	publisher = {Springer},
	title = {Smooth strongly convex interpolation and exact worst-case performance of first-order methods},
	volume = {161},
	year = {2017}}

@article{DT14,
	author = {Drori, Yoel and Teboulle, Marc},
	journal = {Mathematical Programming},
	number = {1},
	pages = {451--482},
	publisher = {Springer},
	title = {Performance of first-order methods for smooth convex minimization: a novel approach},
	volume = {145},
	year = {2014}}

@inproceedings{KNS16,
	author = {Karimi, Hamed and Nutini, Julie and Schmidt, Mark},
	booktitle = {Joint European Conference on Machine Learning and Knowledge Discovery in Databases},
	organization = {Springer},
	pages = {795--811},
	title = {Linear convergence of gradient and proximal-gradient methods under the {Polyak--{\L}ojasiewicz} condition},
	year = {2016}}

@article{P63,
	author = {Polyak, Boris T},
	journal = {Zhurnal Vychislitel'noi Matematiki i Matematicheskoi Fiziki},
	number = {4},
	pages = {643--653},
	title = {Gradient methods for minimizing functionals},
	volume = {3},
	year = {1963}}

@article{THG17b,
	author = {Taylor, Adrien B and Hendrickx, Julien M and Glineur, Fran{\c{c}}ois},
	journal = {SIAM Journal on Optimization},
	number = {3},
	pages = {1283--1313},
	publisher = {SIAM},
	title = {Exact worst-case performance of first-order methods for composite convex optimization},
	volume = {27},
	year = {2017}}

@article{AB24,
	author = {Luner, Alan and Grimmer, Benjamin},
	journal = {arXiv preprint arXiv:2410.14811},
	title = {Performance estimation for smooth and strongly convex sets},
	year = {2024}}

@article{NJF24,
	author = {Bousselmi, Nizar and Hendrickx, Julien M and Glineur, Fran{\c{c}}ois},
	journal = {SIAM Journal on Optimization},
	number = {3},
	pages = {3033--3063},
	publisher = {SIAM},
	title = {Interpolation conditions for linear operators and applications to performance estimation problems},
	volume = {34},
	year = {2024}}

@article{AJ26,
	author = {Rubbens, Anne and Hendrickx, Julien M},
	journal = {Mathematical Programming},
	pages = {1--34},
	publisher = {Springer},
	title = {A constraint-based approach to function interpolation, with application to performance estimation for weakly convex optimization.},
	year = {2026}}

@article{FI18,
	author = {Fikioris, George},
	journal = {Linear Algebra and its Applications},
	pages = {182--210},
	publisher = {Elsevier},
	title = {Spectral properties of {K}ac--{M}urdock--{S}zeg{\"o} matrices with a complex parameter},
	volume = {553},
	year = {2018}}
\newpage
\appendix

\section*{Appendix}

\section{Duality for the primal--dual pair~\eqref{primal:reSDP} and~\eqref{Dual:reSDP}}

\begin{lemma}\label{Lemma:primal_dual}
The Lagrange dual of~\eqref{primal:reSDP} is~\eqref{Dual:reSDP}.
For every primal feasible point $(\ell,\{\delta_k\}_{k=0}^N,\BG)$ and every
 dual feasible point
$(\{\lambda_k,\lambda'_k\}_{k=0}^{N-1},
\{\alpha_k,\beta_k,\eta_k\}_{k=0}^{N},\gamma)$, weak duality gives
\[
\ell\leq \gamma\Delta.
\]
In particular, a primal--dual feasible pair with
$\ell=\gamma\Delta$ is optimal for both programs. Moreover,
~\eqref{primal:reSDP} satisfies Slater's condition; consequently,
~\eqref{primal:reSDP} and~\eqref{Dual:reSDP} have the same optimal value.
\end{lemma}

\begin{proof}
\textbf{Lagrange dual.} We
retain $\BG\succeq 0$ as the conic domain and associate nonnegative
multipliers
\[
\lambda_k,\lambda'_k\geq 0,\quad 0\leq k\leq N-1,
\qquad
\alpha_k,\beta_k,\eta_k\geq 0,\quad 0\leq k\leq N,
\qquad
\gamma\geq 0
\]
with~\eqref{PP1}, \eqref{PP2}, the two inequalities in~\eqref{PP3},
\eqref{PP5}, and~\eqref{PP4}, respectively. Written using nonnegative
primal slacks, the Lagrangian for the maximization problem is
\begin{align}
\mathcal{L}
={}&\ell
+\sum_{k=0}^{N-1}\lambda_k
\left(\delta_k-\delta_{k+1}
      -\lag\BA_{k,k+1},\BG\rag\right)\nonumber\\
&+\sum_{k=0}^{N-1}\lambda'_k
\left(\delta_{k+1}-\delta_k
      -\lag\BA_{k+1,k},\BG\rag\right)\nonumber\\
&+\sum_{k=0}^{N}\alpha_k
\left(\delta_k-\frac{G_{k,k}}{2L}\right)
+\sum_{k=0}^{N}\beta_k
\left(\frac{G_{k,k}}{2\mu}-\delta_k\right)\nonumber\\
&+\gamma(\Delta-\delta_0)
+\sum_{k=0}^{N}\eta_k(\delta_k-\ell).
\label{eq:duality_lagrangian}
\end{align}
Define
\[
\BS:=\sum_{k=0}^{N-1}
\left(\lambda_k\BA_{k,k+1}+\lambda'_k\BA_{k+1,k}\right)
+\diag\left(
\frac{\balpha}{2L}-\frac{\bbeta}{2\mu}
\right).
\]
Collecting the coefficients of $\ell$, $\delta_0,\ldots,\delta_N$, and
$\BG$ in~\eqref{eq:duality_lagrangian} gives
\begin{align}
\mathcal{L}
={}&\gamma\Delta
+\left(1-\sum_{k=0}^{N}\eta_k\right)\ell
-\lag\BS,\BG\rag\nonumber\\
&+\left(
\alpha_0-\beta_0+\eta_0+\lambda_0-\lambda'_0-\gamma
\right)\delta_0\nonumber\\
&+\sum_{k=1}^{N-1}
\left(
\alpha_k-\beta_k+\eta_k-\lambda_{k-1}+\lambda_k
+\lambda'_{k-1}-\lambda'_k
\right)\delta_k\nonumber\\
&+\left(
\alpha_N-\beta_N+\eta_N-\lambda_{N-1}+\lambda'_{N-1}
\right)\delta_N.
\label{eq:duality_lagrangian_collected}
\end{align}
Because $\ell$ and $\delta_0,\ldots,\delta_N$ are unrestricted, and
$\BG$ ranges over the positive semidefinite cone, the supremum of the
Lagrangian is finite if and only if the coefficients of $\ell$ and all
$\delta_k$ vanish and $\BS\succeq 0$. These conditions are exactly
\eqref{DD1}--\eqref{DD3}. Under them,
\[
\sup_{\ell,\delta,\,\BG\succeq 0}\mathcal{L}=\gamma\Delta,
\]
which yields~\eqref{Dual:reSDP}.\\
\\
\textbf{Weak duality and gap identity.} For later use, the same calculation gives the following explicit
primal--dual gap identity for every primal--dual feasible pair:
\begin{align}
\gamma\Delta-\ell
={}&\sum_{k=0}^{N-1}\lambda_k
\left(\delta_k-\delta_{k+1}
      -\lag\BA_{k,k+1},\BG\rag\right)\nonumber\\
&+\sum_{k=0}^{N-1}\lambda'_k
\left(\delta_{k+1}-\delta_k
      -\lag\BA_{k+1,k},\BG\rag\right)\nonumber\\
&+\sum_{k=0}^{N}\alpha_k
\left(\delta_k-\frac{G_{k,k}}{2L}\right)
+\sum_{k=0}^{N}\beta_k
\left(\frac{G_{k,k}}{2\mu}-\delta_k\right)\nonumber\\
&+\gamma(\Delta-\delta_0)
+\sum_{k=0}^{N}\eta_k(\delta_k-\ell)
+\lag\BS,\BG\rag.
\label{eq:primal_dual_gap_identity}
\end{align}
Every term on the right-hand side is nonnegative. In particular,
$\lag\BS,\BG\rag\geq 0$ because $\BS\succeq 0$ and $\BG\succeq 0$.
This proves weak duality and shows directly that matching feasible
objectives certify optimality; no strong-duality argument is needed for
such a certificate.\\
\\
\textbf{Strict primal feasibility.} It remains to verify Slater's condition. Choose $N+1$ distinct numbers
\[
\mu<h_0<h_1<\cdots<h_N<L,
\qquad h_r\neq \frac{1}{t},\quad 0\leq r\leq N,
\]
and set
\[
\BH:=\diag(h_0,\ldots,h_N).
\qquad
0<a<\left(\frac{2\Delta}{\sum_{r=0}^{N}h_r}\right)^{1/2},
\qquad
x_0:=a(1,\ldots,1)^{\top}.
\]
Consider the quadratic function $f(\bx)=\frac12 \bx^{\top}\BH \bx$, whose
minimum value is zero, and generate
\[
\bx_k=(\BI-t\BH)^k \bx_0,
\qquad
\bg_k=\BH \bx_k,
\qquad
\delta_k=\frac12 \bx_k^{\top}\BH \bx_k,
\qquad 0\leq k\leq N.
\]
Let $G_{i,j}:=\lag \bg_i,\bg_j\rag$ and set $\ell:=0$. Since
$\mu\BI\prec\BH\prec L\BI$, for every $k$,
\begin{align*}
\delta_k-\frac{G_{k,k}}{2L}
&=\frac12 \bx_k^{\top}
\left(\BH-\frac{\BH^2}{L}\right)\bx_k>0,\\
\frac{G_{k,k}}{2\mu}-\delta_k
&=\frac12 \bx_k^{\top}
\left(\frac{\BH^2}{\mu}-\BH\right)\bx_k>0.
\end{align*}
Moreover, direct substitution into the definition of $\BA_{i,j}$ gives,
for every ordered pair $i\neq j$,
\begin{equation}
\delta_i-\delta_j-\lag\BA_{i,j},\BG\rag
=\frac{1}{4L}(\bx_i-\bx_j)^{\top}
\left(L^2\BI-\BH^2\right)(\bx_i-\bx_j).
\label{eq:strict_quadratic_interpolation}
\end{equation}
In particular, the right-hand side is strictly positive for the
consecutive pairs appearing in~\eqref{PP1}--\eqref{PP2}, because
$\bx_{k+1}-\bx_k=-t\bg_k\neq 0$ and $L^2\BI-\BH^2\succ 0$.
The choice of $a$ gives $\delta_0<\Delta$, while
$\delta_k>0=\ell$ for all $k$.

Finally, let $\BV=[\bg_0\ \cdots\ \bg_N]$. With
$q_r:=1-th_r$, its entries are
\[
V_{r,k}=a h_r q_r^k,
\qquad 0\leq r,k\leq N.
\]
The numbers $q_0,\ldots,q_N$ are distinct, so $\BV$ is the product of
an invertible diagonal matrix and a Vandermonde matrix. Hence $\BV$ is
invertible and
\[
\BG=\BV^{\top}\BV\succ 0.
\]
Thus every scalar inequality in~\eqref{primal:reSDP} is strict and
$\BG\succ0$, so the primal satisfies Slater's condition. Since
$\ell\leq\delta_0\leq\Delta$, its optimal value is finite. The standard
Slater theorem for semidefinite programs therefore gives equality of the
primal and dual optimal values.
\end{proof}

\begin{remark}
The explicit certificates in Theorems~\ref{thm:primal_dual_equal},
\ref{thm:primal_dual_equal_case2},
\ref{thm:primal_dual_equal_case31}, and
\ref{thm:primal_dual_equal_case32} require only the weak-duality part of
Lemma~\ref{Lemma:primal_dual}: once primal feasibility, dual feasibility,
and $\ell=\gamma\Delta$ are verified, optimality follows immediately
from~\eqref{eq:primal_dual_gap_identity}.
\end{remark}
\section{Lemmas for Theorem~\ref{thm:primal_dual_equal} }
In this section we provide lemmas along with proofs that are needed for proving Theorem~\ref{thm:primal_dual_equal}.  All the notations follow the setting in Theorem~\ref{thm:primal_dual_equal}. Notice in this regime we have
\[
0< Lt \leq 1.
\]
For readability we restate the definition of primal and dual variables in this regime and they will be used in the following lemmas in this section and we assume $\Delta=1$. 
We have primal variables given by
\begin{equation}\label{R1_proof:primal_def}
\delta_k = \Delta p^{k},~0\leq k \leq N,~
\BG = 2\mu \Delta 
\begin{bmatrix}
1 \\
\vdots \\
p^{N/2}
\end{bmatrix}
\begin{bmatrix}
1 \\
\vdots \\
p^{N/2}
\end{bmatrix}^{\top}.
\end{equation}
We have dual variables given by
\begin{equation}\label{R1_proof:dual_def}
\begin{aligned}
\alpha_k & = 0,~~0 \leq k\leq N,   \\
\beta_0 & =p^{N}(rp^{-1}-1),~~\beta_{N} = 1-r, \\
\beta_k & = rp^{N-k}(p^{-1}-1),~~\forall ~1 \leq k\leq N-1, \\
\lambda_k & = rp^{N-k-1}, \lambda'_{k}=0~~\forall ~0\leq k\leq N-1,
\end{aligned}
\end{equation}
and $p$ given by
\begin{equation}\label{R1_proof:def_p}
p = \frac{s^2+2(\kappa-1)}{2\kappa+1}
\left(
\frac{
\sqrt{2(\kappa+1)s^2+2(\kappa-1)(2\kappa+1)}+|s|
}{
\sqrt{2(\kappa+1)s^2+2(\kappa-1)(2\kappa+1)}-|s|
}
\right)
\end{equation}
where
\begin{equation}\label{R1_proof:def_rho}
r:= \frac{2 \kappa (s^2+2\kappa-2) }{2(\kappa+1)s^2 + 2(\kappa-1)(2\kappa+1) - \sqrt{2(\kappa+1)s^2+2(\kappa-1)(2\kappa+1)}|s|}.
\end{equation}
\\
The following lemmas are used in the proof of Theorem~\ref{thm:primal_dual_equal}. Lemma~\ref{R1:L1} shows that the square root of convergence rate $p$ is a root to a quadratic equation and it is smaller than $r$. Lemma~\ref{R1:L3} shows that the certificate matrix $\BS$ is positive semidefinite, which verifies~\eqref{DD3}. Lemma~\ref{R1:L4} verifies that smoothness interpolation constraints hold for all pairs $(i,j),~0\leq i \neq j \leq N$, which shows the solution construction is also feasible to~\eqref{def:sdp_pep}, especially the interpolation inequalities are tight for forward consecutive pairs. 
\begin{lemma}\label{R1:L1}
Consider
\begin{equation}\label{eq:psolqua}
h(x):=(2\kappa+1)x^2 + 2sx +( 2 -s^2 - 2\kappa)=0.
\end{equation}
By setting
\begin{equation}\label{eq:q_def}
q =  \frac{\sqrt{ 2(\kappa+1)s^2 + 2 (\kappa-1)(2\kappa+1)} - s }{2\kappa+1}
\end{equation}
i.e., larger root to~\eqref{eq:psolqua}, 
and then it holds that
\[
p = q^2<1,~~r <1,~~rp^{-1}>1.
\]
\end{lemma}
\begin{proof}
\textbf{Identify $p=q^2$.} Define
\[
z:=\sqrt{2(\kappa+1)s^2+2(\kappa-1)(2\kappa+1)},~~z^2 - s^2 = (2\kappa+1)s^2+2(\kappa-1)(2\kappa+1).
\]
By~\eqref{R1_proof:def_p}, it holds that
\begin{equation}
\begin{aligned}
p& = \frac{s^2+2\kappa-2}{2\kappa+1} \cdot \frac{z + |s|}{z - |s|} = \frac{s^2+2\kappa-2}{2\kappa+1} \cdot \frac{(z + |s|)^2}{(2\kappa+1) (s^2+2\kappa-2) } \\
&=\frac{(z + |s|)^2}{(2\kappa+1)^2}=\left(\frac{\sqrt{2(\kappa+1)s^2 + 2(\kappa-1)(2\kappa+1)}-s}{2\kappa+1}\right)^2 = q^2
\end{aligned}
\end{equation}
where $q$ is defined as one root to~\eqref{eq:psolqua}, $s \leq 0$ in this regime.\\
\\
\textbf{Prove $p<1$.} To show $p<1$, it suffices to verify $0 < q < 1$ for all $-1<s \leq 0$, we notice 
$$
h(0) = 2-s^2-2\kappa <0,~~h(1)=3+2s-s^2 >0,~~h(-1) = 3-2s-s^2 >0.
$$
Therefore by continuity both roots of $h$ lie in $(-1,1)$, hence $0<q<1$ and $p<1$.\\
\\
\textbf{Prove $r<1$ and $r>p$.} We have the following equivalence to show $r < 1$ for $-1 < s \leq 0$:
\begin{align*}
 r < 1 & \Longleftrightarrow 2\kappa (s^2+2\kappa-2) < s^2+(2\kappa+1)(s^2+2\kappa-2) - \sqrt{s^4+(2\kappa+1)(s^2+2\kappa-2)s^2} \\
 & \Longleftrightarrow (2s^2+2\kappa-2)^2 > s^4+(2\kappa+1)(s^2+2\kappa-2)s^2  \\
 & \Longleftrightarrow (2s^2+2\kappa-2)^2-s^4 > (2\kappa+1)(s^2+2\kappa-2)s^2 \\
 & \Longleftrightarrow 3s^2+2\kappa-2 > (2\kappa+1)s^2 \\
 & \Longleftrightarrow 2\kappa-2 > (2\kappa-2)s^2 \\
 & \Longleftrightarrow 1 > s^2.
\end{align*}
To show $rp^{-1}>1$, we notice 
$$
r = \frac{2\kappa (s^2+2\kappa-2)}{z^2-z|s|}.
$$
Hence
$$
rp^{-1} = \frac{2\kappa (s^2+2\kappa-2)}{z^2-z|s|} \cdot \frac{2\kappa+1}{s^2+2\kappa-2} \cdot \frac{z - |s|}{z + |s|} = \frac{2\kappa(2\kappa+1)}{z(z+|s|)} > 1, 
$$
since $|s|<1$, and we have
$$
z=\sqrt{2(\kappa+1)s^2+2(\kappa-1)(2\kappa+1)} < \sqrt{2(\kappa+1)+2(\kappa-1)(2\kappa+1)} = 2\kappa.
$$
\end{proof}

\begin{lemma}\label{R1:L3}
The certificate matrix $\BS \in \mathbb{R}^{(N+1) \times (N+1)}$,
\[
\mathbf S
 :=\sum_{k=0}^{N-1}
 \left(
 \lambda_k\BA_{k,k+1}
 +\lambda_k^{\prime}\BA_{k+1,k}
 \right)
 +\diag\left(
 \frac{\balpha}{2L}-\frac{\bbeta}{2\mu}
 \right)
 \succeq 0,~~\forall N\geq 1,
\]
with the choice of dual variables in~\eqref{R1_proof:dual_def}.
\end{lemma}
\begin{proof}
\textbf{Weighted decomposition.} Fix $N\geq1$ and suppress the superscript $(N)$ throughout the proof.
Let $q:=\sqrt p$, as in Lemma~\ref{R1:L1}.
We prove the following decomposition of $\BS$:
\begin{equation}\label{r1:decompS}
\BS
=\sum_{k=0}^{N-1}\omega_k
  (q\be_k-\be_{k+1})(q\be_k-\be_{k+1})^{\top},
\end{equation}
where
\[
\omega_k:=-\frac{rs}{4Lq}\,p^{N-k-1},
\qquad 0\leq k\leq N-1.
\]
\textbf{Scalar identities.} We first record two scalar identities used in the entrywise
verification. By Lemma~\ref{R1:L1}, $p=q^2$ and
\begin{equation}\label{eq:r1:slack-quadratic}
(2\kappa+1)p+2sq+2-s^2-2\kappa=0.
\end{equation}
Moreover, set
\[
z:=\sqrt{2(\kappa+1)s^2+2(\kappa-1)(2\kappa+1)}
  =(2\kappa+1)q+s.
\]
Since $s\leq0$ and
$z^2-s^2=(2\kappa+1)(s^2+2\kappa-2)$,
the definition~\eqref{R1_proof:def_rho} gives
\begin{align*}
r
&=\frac{2\kappa(s^2+2\kappa-2)}{z(z+s)}
 =\frac{2\kappa q}{z},\\
\frac{2\kappa}{r}
&=\frac{z}{q}
 =2\kappa+1+\frac{s}{q}.
\end{align*}
\textbf{Entrywise verification.} With the multipliers in~\eqref{R1_proof:dual_def}, we have
$\balpha=\bzero$ and $\lambda'_k=0$ for $0\leq k\leq N-1$.
Thus $\BS$ is tridiagonal. We now verify~\eqref{r1:decompS}
by comparing the diagonal and off-diagonal entries.
For the first diagonal entry, we have
\begin{align*}
S_{0,0}
&=\frac{(2-s^2)\lambda_0}{4L}-\frac{\beta_0}{2\mu}\\
&=\frac{rp^{N-1}}{4L}
  \left(2-s^2-2\kappa+\frac{2\kappa p}{r}\right)\\
&=-\frac{rsq}{4L}\,p^{N-1}
 =q^2\omega_0.
\end{align*}
For an interior index $1\leq k\leq N-1$, we obtain
\begin{align*}
S_{k,k}
&=\frac{\lambda_{k-1}+(2-s^2)\lambda_k}{4L}
  -\frac{\beta_k}{2\mu}\\
&=\frac{rp^{N-k-1}}{4L}
  \left[p+2-s^2-2\kappa(1-p)\right]\\
&=-\frac{rsq}{2L}\,p^{N-k-1}
 =\omega_{k-1}+q^2\omega_k,
\end{align*}
where we use~\eqref{eq:r1:slack-quadratic} and
$\omega_{k-1}=p\omega_k=q^2\omega_k$.
The terminal diagonal entry contains the multiplier
$\beta_N=1-r$. Hence
\begin{align*}
S_{N,N}
&=\frac{\lambda_{N-1}}{4L}-\frac{\beta_N}{2\mu}=\frac{r}{4L}-\frac{1-r}{2\mu}=\frac{r}{4L}\left(1+2\kappa-\frac{2\kappa}{r}\right)=-\frac{rs}{4Lq}=\omega_{N-1}.
\end{align*}
For the off-diagonal entries, with $0\leq k\leq N-1$, we have
\begin{align*}
S_{k,k+1}=S_{k+1,k}=\frac{s\lambda_k}{4L}=\frac{rs}{4L}\,p^{N-k-1}=-q\omega_k.
\end{align*}
All remaining entries vanish. These calculations establish
\eqref{r1:decompS}.\\
\\
\textbf{Weight signs.} Since $-1<s\leq0$, $r>0$, $q>0$, and $p>0$, every weight
$\omega_k$ is nonnegative. Consequently, for every
$\bx=(x_0,\ldots,x_N)^{\top}\in\mathbb R^{N+1}$,
\[
\bx^{\top}\BS\bx
=\sum_{k=0}^{N-1}\omega_k(qx_k-x_{k+1})^2
\geq0.
\]
Therefore $\BS\succeq0$ for every $N\geq1$.
In particular, at $t=1/L$ we have $s=0$, so every weight
vanishes and $\BS=0$.
\end{proof}

\begin{lemma}\label{R1:L4}
We have for any pair $(i,j),~(0 \leq i \neq j\leq N)$, it holds that
\begin{equation}
\Tr(\BA_{i,j}\BG) \leq \delta_i - \delta_j.
\end{equation}
\end{lemma}
\begin{proof}
\textbf{General rank-one trace identities.} We start by recording identities that also apply to the rank-one construction in Regime~III. This calculation only concerns the inner product between $\BA_{i,j}$ and a symmetric matrix $\BM$ with a rank one geometric decay structure.
Fix $0\leq i<j\leq N$, put $h:=j-i$, and let a symmetric matrix
$\BM$ satisfy
\[
 M_{i+r,\,i+u}=H\rho^{r+u},
 \qquad 0\leq r,u\leq h,
\]
where $H>0$ and $-1<\rho<1$. 
Finite geometric sums give
\begin{equation}\label{eq:rankone:geometric-sums}
\begin{aligned}
 M_{i,i}+M_{j,j}-2M_{i,j}
 &=H(1-\rho^h)^2,\\
 \sum_{r=i}^{j-1}(M_{i,r}+M_{j,r})
 &=H(1+\rho^{h})\sum_{a=0}^{h-1}\rho^a
 =\frac{H(1-\rho^{2h})}{1-\rho},\\
 \sum_{r=i}^{j-1}\sum_{u=i}^{j-1}M_{r,u}
 &=H\left(\sum_{a=0}^{h-1}\rho^a\right)^2
 =\frac{H(1-\rho^{h})^2}{(1-\rho)^2}.
\end{aligned}
\end{equation}
The direct definition of $\BA_{i,j}$ \eqref{def:Aij} now yields
\begin{equation}\label{eq:rankone:pair-traces}
\begin{aligned}
 \lag\BA_{i,j},\BM\rag
 &=\frac{H}{4L}\left(1-\frac{L^2t^2}{(1-\rho)^2}\right)
      (1-\rho^{h})^2
   +\frac{tH}{2(1-\rho)}(1-\rho^{2h}),\\
 \lag\BA_{j,i},\BM\rag
 &=\frac{H}{4L}\left(1-\frac{L^2t^2}{(1-\rho)^2}\right)
      (1-\rho^{h})^2
   -\frac{tH}{2(1-\rho)}(1-\rho^{2h}).
\end{aligned}
\end{equation}
Only the term linear in $t$ changes sign when the ordered
pair is reversed. Identities~\eqref{eq:rankone:geometric-sums} and \eqref{eq:rankone:pair-traces} remain valid for negative $\rho$.\\
\\
\textbf{Regime-I specialization and scalar signs.} We now specialize the preceding calculation to the
Regime~I construction. Recall $q:=\sqrt p$ from~\eqref{eq:q_def}. Then
\[
G_{k,l}=2\mu q^{k+l},
\qquad
\delta_k= q^{2k},
\qquad
0\leq k,l\leq N.
\]
By Lemma~\ref{R1:L1}, we have $0<q<1$, and $q$ is the
positive root of
\[
(2\kappa+1)q^2+2sq+2-s^2-2\kappa=0.
\]
We first record two consequences of this quadratic identity.
Since $0<\mu<L$ and $0<t\leq1/L$, we have $1-\mu t>0$.
Evaluating the polynomial in~\eqref{eq:psolqua} at
$1-\mu t$ gives
\begin{equation}\label{eq:h_mut}
\begin{aligned}
(2\kappa+1)(1-\mu t)^2   +2s(1-\mu t)+2-s^2-2\kappa =(\mu^2-L^2)t^2<0,
\end{aligned}
\end{equation}
where $s=Lt-1$ and $\kappa=L/\mu$.
Because this polynomial has a positive leading coefficient
and exactly one positive root, it follows that
\[
q>1-\mu t.
\]
Therefore, we have
\[
0<1-q<\mu t<Lt
\]
and 
\begin{equation}\label{eq:r1:positive-chi}
\chi := \frac{L^2t^2-(1-q)^2}{4L(1-q)^2} >0.
\end{equation}
Moreover, substituting $s=Lt-1$ into the quadratic
identity for $q$ and rearranging yields
\[
2\kappa(1-q^2)
=
(1-q)^2-L^2t^2+2Lt(1+q).
\]
Dividing by $4L(1-q^2)$, we obtain
\begin{equation}\label{eq:r1:slack-coefficient}
\begin{aligned}
\frac{1}{2\mu}-\frac{t}{2(1-q)}
&=
\frac{(1-q)^2-L^2t^2}{4L(1-q^2)}\\
&=-\chi\,\frac{1-q}{1+q}.
\end{aligned}
\end{equation}
Fix $0\leq i<j\leq N$, and set $h:=j-i$ and
$H:=G_{i,i}>0$. The entries of $\BG$ satisfy
\[
G_{i+r,i+u}=Hq^{r+u},
\qquad 0\leq r,u\leq h,
\]
while
\[
\delta_i-\delta_j
=
\frac{H}{2\mu}(1-q^{2h}).
\]
Applying~\eqref{eq:rankone:pair-traces} with
$\BM=\BG$ and $\rho=q$ gives
\begin{align*}
\lag\BA_{i,j},\BG\rag
&=-\chi H(1-q^h)^2
  +\frac{tH}{2(1-q)}(1-q^{2h}),\\
\lag\BA_{j,i},\BG\rag
&=-\chi H(1-q^h)^2
  -\frac{tH}{2(1-q)}(1-q^{2h}).
\end{align*}
\textbf{Forward pairs.} For the ordered pair $(i,j)$, using
\eqref{eq:r1:slack-coefficient}, we obtain
\begin{equation}\label{eq:r1:forward-slack}
\begin{aligned}
&\delta_i-\delta_j-\lag\BA_{i,j},\BG\rag\\
&\quad=
\chi H(1-q^h)^2
+
\left(\frac{1}{2\mu}-\frac{t}{2(1-q)}\right)
H(1-q^{2h})\\
&\quad=
\chi H(1-q^h)
\left[
(1-q^h)-\frac{1-q}{1+q}(1+q^h)
\right]\\
&\quad=
\frac{2\chi H}{1+q}(1-q^h)(q-q^h)
\geq0.
\end{aligned}
\end{equation}
Indeed, $0<q<1$ and $h\geq1$ imply $q^h\leq q$.
Equality holds precisely when $h=1$, whereas the
inequality is strict when $h\geq2$.\\
\\
\textbf{Reversed pairs.} For the reversed ordered pair $(j,i)$, the same
coefficient identity yields
\begin{equation}\label{eq:r1:reverse-slack}
\begin{aligned}
&\delta_j-\delta_i-\lag\BA_{j,i},\BG\rag\\
&\quad=
\chi H(1-q^h)^2
+
\left(\frac{t}{2(1-q)}-\frac{1}{2\mu}\right)
H(1-q^{2h})\\
&\quad=
\chi H(1-q^h)
\left[
(1-q^h)+\frac{1-q}{1+q}(1+q^h)
\right]\\
&\quad=
\frac{2\chi H}{1+q}(1-q^h)(1-q^{h+1})
>0.
\end{aligned}
\end{equation}
Thus, the interpolation constraints hold for every
ordered pair. In particular, the forward consecutive
constraints hold with equality, and all other
ordered-pair constraints hold strictly.
\end{proof}

\section{Lemmas for Theorem~\ref{thm:primal_dual_equal_case2} }
In this section we provide lemmas along with proofs that are needed for proving Theorem~\ref{thm:primal_dual_equal_case2}.  All the notations follow the setting in Theorem~\ref{thm:primal_dual_equal_case2}. Hence, $\delta_k,\BG$ are defined in~\eqref{R2_proof:vari_primal_def}, $\alpha_k,\beta_k,\lambda_k$ are defined in~\eqref{R2_proof:vari_dual_def}. Notice in this regime we have
$$
1< Lt \leq \sqrt{3}.
$$
For readability we restate the definitions of primal and dual variables in this regime and they will be used in the following lemmas in this section and we assume $\Delta=1$. For the primals we have
\begin{equation}\label{R2_proof:vari_primal_def}
\begin{aligned}
\delta_k = p^{k}\Delta,~\forall~0\leq k\leq N,~\ell = p^{N}\Delta, ~~~~ G^{(N)} = 2\mu\Delta \BP^{(N)} \BPhi^{(N)} \BP^{(N)},
\end{aligned}
\end{equation}
where
\begin{equation}\label{def_proof:BPBPHI}
\begin{aligned}
\BP^{(N)}
=\diag\left(1,p^{1/2},\ldots,p^{N/2}\right),
\qquad
\BPhi^{(N)}
=\left(\phi_{|i-j|}\right)_{i,j=0}^{N},
\end{aligned}
\end{equation}
and $\phi_k$'s are given by the following recursive relation
\begin{equation}\label{R2_proof:eq_recur}
\phi_0=1,~~\phi_1=\frac{\tau}{\sqrt{p}},~~\phi_{k+2} = \frac{s}{s+2}\phi_k + \frac{p-s(s+2)}{\sqrt{p}(s+2)}\phi_{k+1},~k\geq 0
\end{equation}
where
\begin{equation}\label{R2_proof:p_def}
p := 1+\frac{(s^2-1)(s+3)}{s+2\kappa+1} = \frac{\mu L^2t^3-3\mu t+2}{2 + \mu t} ,~~\tau=\frac{1}{2}(1-L^2t^2+p) = \frac{2-L^2t^2-\mu t}{2+\mu t}.
\end{equation}
For the duals we have
\begin{equation}\label{R2_proof:vari_dual_def}
\begin{aligned}
&\lambda_{k} = p^{N-k-1} \frac{s+2}{2+\mu t},~~ \lambda_k' = p^{N-k-1} \frac{s}{2+\mu t},~~0\leq k\leq N-1, \\
&\alpha_{k} = 0,~~0\leq k\leq N,\\
&\beta_k = \lambda'_{k-1}-\lambda'_k - (\lambda_{k-1}-\lambda_k),~1\leq k\leq N-1,\\
&\beta_0 = \lambda_0-\lambda_0'-p^N,~~\beta_{N} = 1 + \lambda_{N-1}' - \lambda_{N-1}. \\
\end{aligned}
\end{equation}
The following lemmas are used in the proof of Theorem~\ref{thm:primal_dual_equal_case2}. Lemma~\ref{R2:L1} shows that $0<p<1$. Lemma~\ref{R2:L2} shows that the primal solution $\BG$ satisfies the constraints \eqref{PP1} for all pairs, hence the optimal solution is also feasible for the full SDP~\eqref{def:sdp_pep}. 
Lemma~\ref{lemma:I_induct} and~\ref{lemma:recur_cf} facilitate the proof of Lemma~\ref{R2:L2}. Lemma~\ref{R2:L3} and~\ref{R2:L4} establish positive semidefiniteness of $\BG$ and $\BS$ resepctively. \\
\\
In particular, the primal-feasibility results in
Lemmas~\ref{R2:L1}, \ref{R2:L2}, and~\ref{R2:L3},
together with the supporting scalar results in
Lemmas~\ref{lemma:I_induct} and~\ref{lemma:recur_cf},
hold throughout $1/L<t<2/L$, not only in Regime~II.
They are therefore also used in
Proposition~\ref{prop:r3c2:toeplitz-full-feasible}.
\begin{lemma}\label{R2:L1}
We have $s^2<p<1$ for $1/L<t<2/L$.
\end{lemma}
\begin{proof}
On the one hand we have
\begin{align*}
p-s^2 &= \frac{\mu L^2t^3-3\mu t+2}{2 + \mu t} - (1-Lt)^2 \\
& = \frac{\mu L^2t^3-3\mu t+2-(2+\mu t)(Lt-1)^2}{2 + \mu t} \\
& = \frac{2t(L-\mu)(2-Lt)}{2+\mu t} > 0,
\end{align*}
since $Lt<2$ and $L>\mu$. On the other hand we have
$$
1-p = \frac{\mu t (4-L^2t^2)}{2+\mu t} > 0,
$$
since $1<Lt<2$.
\end{proof}

\begin{lemma}\label{R2:L2}
Fix $N>0$ and $1/L<t<2/L$, we have for $(i,j)$ pair with $0 \leq i,j \leq N$,
\begin{enumerate}[label=(\alph*)]
\item \label{R2:L2S1}if $|i-j| = 1$, then $\Tr(\BA_{i,j}\BG) = \delta_i - \delta_j$,

\item \label{R2:L2S2}if $i>j$, then $\Tr(\BA_{i,j}\BG)=\delta_{i}-\delta_{j}$,

\item \label{R2:L2S3}if $j\geq i+2$, then $\Tr(\BA_{i,j}\BG)<\delta_{i}-\delta_{j}$.
\end{enumerate}
We prove these statements one by one.
\end{lemma}
\begin{proof} \textbf{Proof for \ref{R2:L2S1}}: Notice for fixed $k$ we have
\begin{align*}
\Tr(\BA_{k,k+1}\BG) = p^{k} \Tr(\BA_{0,1}\BG),~~p^{k}(\delta_0-\delta_1) = \delta_{k}-\delta_{k+1}, \\
\Tr(\BA_{k+1,k}\BG) = p^{k} \Tr(\BA_{1,0}\BG), ~~p^{k}(\delta_1-\delta_0) = \delta_{k+1}-\delta_{k}.
\end{align*}  
Hence it suffices to prove for the statement $(0,1)$ and $(1,0)$ pair, which is equivalent to have
\begin{align*}
\Tr((\BA_{0,1}+\BA_{1,0})\BG) &= 0, \\
\Tr((\BA_{0,1}-\BA_{1,0})\BG) &= 2(\delta_0-\delta_1) = \frac{1}{\mu}(G_{00}-G_{11}).
\end{align*}
As
\begin{align*}
\Tr((\BA_{0,1}+\BA_{1,0})\BG) &= \Tr\left((\BA+\BA')\BG_{0:1,0:1} \right),~~\Tr((\BA_{0,1}-\BA_{1,0})\BG) = \Tr\left((\BA-\BA')\BG_{0:1,0:1} \right), \\
\BA + \BA' &= \frac{1}{2L}\begin{bmatrix} -s^2-2s & -1 \\ -1 & 1  \end{bmatrix},~~\BA - \BA' = \frac{1}{2L}\begin{bmatrix} 2s+2 & s+1 \\ s+1 & 0  \end{bmatrix}.
\end{align*}
We need the following system of equations hold,
\begin{align*}
\Tr((\BA_{0,1}+\BA_{1,0})\BG) = 0, ~~~~~\Tr((\BA_{0,1}-\BA_{1,0})\BG) = 2(\delta_0-\delta_1).
\end{align*}
This system holds with relation
$$
G_{11} = pG_{00}, G_{01}=G_{10} = \left( \frac{1-p}{\mu t} - 1 \right)G_{00} = \frac{2-L^2t^2-\mu t}{2+\mu t}G_{00}=\tau G_{00} = 2\mu \sqrt{p}\phi_1,
$$
which is consistent with the definition for $\BG$ in~\eqref{R2_proof:vari_primal_def}.\\
\\
\textbf{Proof for \ref{R2:L2S2}}: We prove by induction on $N$. \\
\\
\textbf{Base case:} When $N=1$, the statement holds since we need to prove only pair $(1,0)$, which is covered in~\ref{R2:L2S1}.\\
\\
\noindent\textbf{Induction Step:}
Notice in the way we define $\BG$, suppose we have the statement hold for $N \geq 1$, we want to prove the statement for $N+1$. Notice we have
\begin{align*}
&\BG^{(N+1)}_{0:N,0:N} = \BG^{(N)}, \delta^{(N+1)}_i = \delta^{(N)}_i,~0\leq i \leq N, \\
&\BG^{(N+1)}_{1:N+1,1:N+1} = p\BG^{(N)},\delta^{(N+1)}_{i+1} = p\delta^{(N)}_i,~0\leq i \leq N.
\end{align*}
Hence by induction hypothesis, we have for $(i,j)$ pair with $0\leq j<i \leq N$:
$$
\Tr(\BA^{(N+1)}_{i,j}\BG^{(N+1)}) = \Tr(\BA^{(N)}_{i,j}\BG^{(N)}) = \delta^{(N)}_i - \delta^{(N)}_j = \delta^{(N+1)}_{i}-\delta^{(N+1)}_j.
$$
For $(i,j)$ pair with $1 \leq j < i \leq N+1$:
$$
\Tr(\BA^{(N+1)}_{i,j}\BG^{(N+1)}) = p\Tr(\BA^{(N)}_{i-1,j-1}\BG^{(N)}) = p\delta^{(N)}_{i-1} - p\delta^{(N)}_{j-1} = \delta^{(N+1)}_{i}-\delta^{(N+1)}_j.
$$
Therefore it suffices to verify only for pair $(N+1,0)$, i.e. we need
$$
\Tr(\BA^{(N+1)}_{N+1,0}\BG^{(N+1)}) = \delta^{(N+1)}_{N+1}-\delta^{(N+1)}_0.
$$
Then it is equivalent to have
\begin{equation}\label{eq:R2_proof_induction}
\begin{aligned}
&\Tr(\BA^{(N+1)}_{N+1,0}\BG^{(N+1)}) - \Tr(\BA^{(N+1)}_{N,0}\BG^{(N+1)}) \\
& = \delta^{(N+1)}_{N+1}-\delta^{(N+1)}_N \\
& =p(\delta^{(N+1)}_{N}-\delta^{(N+1)}_{N-1}) \\
& = p\left[\Tr(\BA^{(N+1)}_{N,0}\BG^{(N+1)}) - \Tr(\BA^{(N+1)}_{N-1,0}\BG^{(N+1)})\right]
\end{aligned}
\end{equation}
where we use convention $\Tr(\BA^{(N+1)}_{0,0}\BG^{(N+1)})=0$. Suppressing the superscript $(N+1)$ for clarity and plugging in equation~\eqref{eq:Aij-direct-inner-product}, 
\begin{align*}
&\frac{G_{N+1,N+1} - G_{N,N} - 2G_{N+1,0} +2G_{N,0}}{4L} + \frac{t}{2}(1-p)\sum_{r=0}^{N-1}G_{N,r} - \frac{t}{2}(G_{N+1,0}+G_{N,0}) \\
& - \frac{Lt^2}{2}\sum_{r=0}^{N-1}G_{N,r} - \frac{Lt^2}{4}G_{N,N} \\
&= ~~\frac{p(G_{N,N} - G_{N-1,N-1} - 2G_{N,0} +2G_{N-1,0})}{4L} + \frac{t}{2}(1-p)p \sum_{r=0}^{N-2} G_{N-1,r} - \frac{tp}{2}(G_{N,0}+G_{N-1,0}) \\
&~~~~ -\frac{Lt^2p}{2}\sum_{r=0}^{N-2} G_{N-1,r} - \frac{Lt^2p}{4}G_{N-1,N-1}.
\end{align*}
Using $G_{i,j}=pG_{i-1,j-1}$ for $i,j\geq1$, we have
\begin{align*}
&G_{N+1,N+1} - G_{N,N} - 2G_{N+1,0} +2G_{N,0} - p(G_{N,N} - G_{N-1,N-1} - 2G_{N,0} +2G_{N-1,0})  \\
=&2((1+p)G_{N,0}-G_{N+1,0}-pG_{N-1,0}), \\
& \sum_{r=0}^{N-1}G_{N,r} - p \sum_{r=0}^{N-2}G_{N-1,r} = \sum_{r=0}^{N-1}G_{N,r} - \sum_{r=1}^{N-1}G_{N,r} = G_{N,0}, \\
& G_{N,N} - pG_{N-1,N-1} = 0.
\end{align*}
Plugging the equalities above into~\eqref{eq:R2_proof_induction} gives the following equivalent equation only containing $G_{N-1,0},G_{N,0},G_{N+1,0}$,
\[
\frac{(1+p)G_{N,0}-G_{N+1,0}-pG_{N-1,0}}{2L} + \frac{t(1-p)}{2}G_{N,0}-\frac{Lt^2}{2}G_{N,0} = \frac{t}{2}(G_{N+1,0}+(1-p)G_{N,0}-pG_{N-1,0}).
\]
Rearranging coefficients gives 
\begin{equation}
\begin{aligned}
&\frac{1}{L}\left(Lt+1\right)G_{N+1,0} = \frac{1}{L}\left(p + 1 - L^2t^2 \right) G_{N,0} + \frac{p(Lt-1)}{L}G_{N-1,0} \\
\Longleftrightarrow&\frac{1}{L}\left(s+2\right)G_{N+1,0} = \frac{1}{L}\left(p - s(s+2) \right) G_{N,0} + \frac{ps}{L}G_{N-1,0}, \\
\Longleftrightarrow &G_{N+1,0} =  \frac{s}{s+2}\cdot p G_{N-1,0} + \frac{p-s(s+2)}{\sqrt{p}(s+2)} \cdot \sqrt{p} G_{N,0}, \\
\Longleftrightarrow &\phi_{N+1} =  \frac{s}{s+2}\cdot \phi_{N-1} + \frac{p-s(s+2)}{\sqrt{p}(s+2)}  \phi_{N},
\end{aligned}
\end{equation}
where we use the fact that $G_{i,0} = 2\mu \sqrt{p}^{i}\phi_i$. The last line is exactly the recursive relation in~\eqref{R2_proof:eq_recur}. Hence the statement holds.\\
\\
\textbf{Proof for \ref{R2:L2S3}}: Following the same induction argument in (2), it suffices show for fixed $N \geq 2$
\begin{equation}
\Tr(\BA_{0,N}\BG) \leq \delta_{0} - \delta_{N}.
\end{equation}
Notice we have
\begin{align*}
&\Tr(\BA_{N,0}\BG) = \delta_{N} - \delta_{0}, \\
&\Tr(\BA_{0,N}\BG) - \Tr(\BA_{N,0}\BG) = t \sum_{k=0}^{N-1}(G_{N,k}+G_{0,k}).
\end{align*}
Using
\[
G_{i,j}
=2\mu p^{(i+j)/2}\phi_{|i-j|},
\qquad
\delta_k= p^k,
\]
it is therefore equivalent to show
\begin{equation}\label{eq:r2_induct1}
2\mu t\sum_{k=0}^{N-1}
\left(
p^{(N+k)/2}\phi_{N-k}
+p^{k/2}\phi_k
\right)
\leq 2(1-p^N).
\end{equation}
We prove this inequality for every $N\geq0$ and show that
it is strict for $N\geq2$.
All sums over empty index sets are understood to be zero.\\
\\
\noindent\textbf{Scalar recurrence.}
Define
\[
H_n
:=
2\mu t\sum_{k=0}^{n-1}
\left(
p^{(n+k)/2}\phi_{n-k}
+p^{k/2}\phi_k
\right),
\qquad n\geq0.
\]
In particular, $H_0=0$.
Reindexing the first sum gives
\[
\begin{aligned}
&\sum_{k=0}^{n}
p^{(n+1+k)/2}\phi_{n+1-k}
-p\sum_{k=0}^{n-1}
p^{(n+k)/2}\phi_{n-k}\\
&\qquad=p^{(n+1)/2}\phi_{n+1}.
\end{aligned}
\]
Therefore, it holds that
\begin{align*}
H_{n+1}-pH_n & =
2\mu t\left[
p^{n/2}\phi_n
+p^{(n+1)/2}\phi_{n+1}
+(1-p)\sum_{k=0}^{n-1}p^{k/2}\phi_k
\right]\\
&=I(n),
\end{align*}
where $I(n)$ is the scalar sequence defined in
Lemma~\ref{lemma:I_induct}.

\medskip
\noindent\textbf{Bound on the increment.}
By Lemma~\ref{lemma:I_induct},
\[
I(n)\leq I(0),
\qquad n\geq0,
\]
with strict inequality for every $n\geq1$.
The initial conditions $\phi_0=1$ and
$\sqrt p\,\phi_1=\tau$ give
\[
I(0)
=2\mu t\left(\phi_0+\sqrt p\,\phi_1\right)
=2\mu t(1+\tau)
=2(1-p),
\]
where the last equality follows from
$\mu t(1+\tau)=1-p$.
Thus,
\[
I(n)\leq2(1-p),
\qquad n\geq0,
\]
with equality only when $n=0$.

\medskip
\noindent\textbf{Induction and strictness.}
The inequality $H_n\leq2(1-p^n)$ holds at $n=0$.
Suppose it holds for some $n\geq0$.
Since $p>0$, the scalar recurrence yields
\[
H_{n+1}
=pH_n+I(n)
\leq2p(1-p^n)+2(1-p)
=2(1-p^{n+1}).
\]
This proves~\eqref{eq:r2_induct1} for all $N\geq0$.

Moreover, $H_1=I(0)=2(1-p)$.
For every $n\geq1$, the strict inequality
$I(n)<2(1-p)$ gives
\[
H_{n+1}
=pH_n+I(n)
<2p(1-p^n)+2(1-p)
=2(1-p^{n+1}).
\]
Hence~\eqref{eq:r2_induct1} is strict for every $N\geq2$.
By the index-shift identities established in
part~\ref{R2:L2S2}, this proves
part~\ref{R2:L2S3} for every $j\geq i+2$.
\end{proof}

\begin{lemma}\label{lemma:I_induct}
Let $1/L<t<2/L$, and let $\{\phi_n\}_{n\geq0}$ be the
sequence defined by~\eqref{R2_proof:eq_recur}.
For every integer $n\geq0$, define
\begin{equation}
I(n)
:=
2\mu t\left[
p^{n/2}\phi_n
+p^{(n+1)/2}\phi_{n+1}
+(1-p)\sum_{k=0}^{n-1}p^{k/2}\phi_k
\right],
\end{equation}
where the sum is zero when $n=0$.
Then, for every even integer $n\geq0$,
\begin{enumerate}[label=(\alph*)]
    \item \label{state:I_induct_1}
    $I(n+2)<I(n)$;
    \item \label{state:I_induct_2}
    $I(n+1)<I(n)$.
\end{enumerate}
Thus, we have
\[
I(n)<I(0),
\qquad n\geq1.
\]
\end{lemma}

\begin{proof}
By Lemma~\ref{lemma:recur_cf}, the sequence admits the
representation
\[
\phi_n=C_1q_1^n+C_2q_2^n,
\qquad n\geq0,
\]
where
\[
C_1,C_2>0,\qquad
q_1=-\frac{s}{\sqrt p}\in(-1,0),
\qquad
q_2=\frac{\sqrt p}{s+2}\in(0,1).
\]
Fix an even integer $n\geq0$.

\medskip
\noindent\textbf{Two-step decrease.}
From the definition of $I(n)$, we obtain
\begin{align*}
I(n+2)-I(n)
&=
2\mu t\left[
p^{(n+3)/2}\phi_{n+3}
+p^{(n+2)/2}\phi_{n+2}
-p^{(n+3)/2}\phi_{n+1}
-p^{(n+2)/2}\phi_n
\right]\\
&=
2\mu t\,p^{(n+2)/2}
\left[
\sqrt p\,(\phi_{n+3}-\phi_{n+1})
+\phi_{n+2}-\phi_n
\right].
\end{align*}
Substituting the closed-form representation gives
\begin{align*}
&\sqrt p\,(\phi_{n+3}-\phi_{n+1})
+\phi_{n+2}-\phi_n\\
&\quad=
C_1q_1^n(1+\sqrt p\,q_1)(q_1^2-1)
+C_2q_2^n(1+\sqrt p\,q_2)(q_2^2-1).
\end{align*}
Since $n$ is even, both $q_1^n$ and $q_2^n$ are positive.
Moreover,
\[
1+\sqrt p\,q_1=1-s>0,
\qquad
1+\sqrt p\,q_2=1+\frac{p}{s+2}>0,
\]
because $0<s<1$.
Since $q_1^2<1$ and $q_2^2<1$, both summands are strictly
negative. Therefore,
\begin{equation}\label{eq:I_inductioneq1}
\sqrt p\,\phi_{n+3}+\phi_{n+2}
<
\sqrt p\,\phi_{n+1}+\phi_n.
\end{equation}
As $2\mu t\,p^{(n+2)/2}>0$, this proves
$I(n+2)<I(n)$, establishing
part~\ref{state:I_induct_1}.

\medskip
\noindent\textbf{Comparison with the next odd index.}
Similarly,
\begin{align*}
I(n+1)-I(n)
&=
2\mu t\,p^{(n+2)/2}
\left(\phi_{n+2}-\phi_n\right)\\
&=
2\mu t\,p^{(n+2)/2}
\left[
C_1q_1^n(q_1^2-1)
+C_2q_2^n(q_2^2-1)
\right].
\end{align*}
Again, $n$ is even, $C_1,C_2>0$, and $|q_1|,|q_2|<1$,
so both summands are strictly negative.
Hence $I(n+1)<I(n)$, proving
part~\ref{state:I_induct_2}.

\medskip
\noindent\textbf{Comparison with the initial value.}
Part~\ref{state:I_induct_1} shows that the even subsequence
is strictly decreasing:
\[
I(2m)<I(2m-2)<\cdots<I(0),
\qquad m\geq1.
\]
For every $m\geq0$, part~\ref{state:I_induct_2} gives
\[
I(2m+1)<I(2m)\leq I(0).
\]
Combining the even and odd cases yields
$I(n)<I(0)$ for every $n\geq1$.
\end{proof}

\begin{lemma}\label{lemma:recur_cf} For recurrence system~\eqref{R2:eq_recur}, the sequence $\{\phi_k\}_{k \geq 0}$ has the closed-form representation
\begin{equation}\label{eq:recur_cf}
\phi_k = C_1 q_1^{k} + C_2q_2^{k}
\end{equation}
where we have $C_1,C_2>0$,$-1<q_1<0$ and $0<q_2<1$. 
\end{lemma}

\begin{proof}
\textbf{Solving the recurrence.} The recurrence system is linear and homogeneous with constant coefficients, the characteristic equation is 
$$
q^2 - \frac{p-s(s+2)}{\sqrt{p}(s+2)}q-\frac{s}{s+2}=0
$$
with roots
$$
q_1 = -\frac{s}{\sqrt{p}},~~q_2 = \frac{\sqrt{p}}{s+2}.
$$
The initial conditions determine
$$
\phi_{k} = \frac{1}{p+s(s+2)} \left[(p-(s+2)\tau) \left(\frac{-s}{\sqrt{p}}\right)^{k}+(\tau+Lt\tau+s(s+2)) \left(\frac{\sqrt{p}}{s+2}\right)^{k}\right].
$$
Hence we have
$$
C_1 =  \frac{p-(s+2)\tau}{p+s(s+2)},~C_2 = \frac{(s+2)(\tau+s)}{p+s(s+2)},~q_1 = -\frac{s}{\sqrt{p}},~q_2=\frac{\sqrt{p}}{s+2}.
$$
\textbf{Positivity of coefficients.} To show $C_1>0,C_2>0$, since $s(s+2)>0$ and $p>0$, we have the denominator is greater than $0$, so it suffices to verify the numerators are greater than $0$. For $C_1$, we have
\begin{align*}
p-(s+2)\tau &= \frac{1}{2+\mu t}\left(\mu L^2t^3+2-3\mu t - (1+Lt)(2-L^2t^2-\mu t) \right) \\
&=\frac{1}{2+\mu t}\left(\mu L^2t^3 +L^2t^2  +L^3t^3+ \mu Lt^2- 2\mu t - 2Lt \right) \\
&> \frac{1}{2+\mu t}\left(\mu t +Lt  +Lt+ \mu t- 2\mu t - 2Lt \right) = 0
\end{align*}
where for the inequality we apply the fact that $Lt>1$. For $C_2$, we have
\begin{align*}
\tau+s  &= \frac{1}{2+\mu t}\left(2Lt+\mu Lt^2 - L^2t^2- 2\mu t \right) \\
&=\frac{1}{2+\mu t}t(L-\mu)(2-Lt) >0
\end{align*}
where for the last line, we apply $L>\mu$ and $Lt<2$. \\
\\
\textbf{Root intervals.} For $q_1$, we have $q_1<0$, since $s>0$, to show $q_1>-1$, it suffices to verify $s < \sqrt{p}$, which is implied by Lemma~\ref{R2:L1}.\\
For $q_2$, we have $q_2>0$ since $s+2=Lt+1>0$ and $q_2<1$ since $\sqrt{p}<1<Lt+1=s+2$ by Lemma~\ref{R2:L1}.
\end{proof}
\begin{lemma}\label{R2:L3}
We have for any $N \geq 1, 1/L < t <2/L, \BG\succ 0$.
\end{lemma}
\begin{proof}
We define the Kac–Murdock–Szegő matrix as
$$
\BK^{(N)}(q):= \left(q^{|i-j|} \right)_{i,j=0}^{N}.
$$
Notice $\BK^{(N)}(q)$ is positive definite for $-1<q<1$, see Theorem 2.3 (v) in \cite{FI18} for reference. Since we have
$$
\BG = 2\mu \BP \BPhi \BP= 2\mu \diag\left(1,p^{1/2},\ldots,p^{N/2}\right) \BPhi \diag\left(1,p^{1/2},\ldots,p^{N/2}\right).
$$
By Lemma~\ref{lemma:recur_cf}, we have 
$$
\BPhi^{(N)} = C_1 \BK^{(N)}(q_1) + C_2 \BK^{(N)}(q_2).
$$
Notice $ -1<q_1<0$ and $0<q_2<1$, hence $\BK^{(N)}(q_1),\BK^{(N)}(q_2) \succ 0$. Since $C_1,C_2>0$, we have $\BPhi \succ 0$. Since $p>0$, the diagonal matrix $\BP^{(N)}$ is invertible.
Therefore, the following holds
\[
\BG = 2\mu\BP^{(N)}\BPhi^{(N)}\BP^{(N)}
\succ0.
\]
\end{proof}
\begin{lemma}\label{R2:L4}
The certificate matrix $\BS \in \mathbb{R}^{N+1 \times N+1}$
$$
\BS = \sum_{k=0}^{N-1} \left(\lambda_k\BA_{k,k+1} + \lambda'_k \BA_{k+1,k}\right) + \diag \left(\frac{\balpha}{2L}-\frac{\bbeta}{2\mu}\right) = 0
$$
with choice of dual variables given by~\eqref{R2_proof:vari_dual_def}.
\end{lemma}
\begin{proof}
Notice $\BS$ is a tri-diagonal matrix. It suffices to verify entries on the tridiagonal band.\\
\\
\textbf{Off-diagonal Entries.} By directly plugging in we have
\[
S_{k,k+1} = S_{k+1,k}= \frac{\lambda_k}{4L}s - \frac{\lambda'_k}{4L}\left(s+2\right) = 0
\]
since $s \lambda_k = (s+2)\lambda_k'$ by~\eqref{R2_proof:vari_dual_def}. \\
\\
\textbf{Diagonal Endpoint Entries.} Now we check elements on the diagonal, on the two endpoints, we have
\begin{align*}
S_{0,0}:&=\lambda_0A_{0,0}+\lambda_0'A'_{0,0}-\frac{1}{2\mu}\beta_0+\frac{1}{2L}\alpha_0\\
&=\frac{\lambda_0}{4L}\left(2-s^2\right)+\frac{\lambda'_0}{4L}\left(-s^2-4s-2\right)-\frac{1}{2\mu}\left( \lambda_0-\lambda_0'-p^N\right), \\
&= \frac{p^{N-1}}{4L} \left[ \frac{s+2}{2+\mu t}\left(2-s^2\right)+ \frac{s}{2+\mu t}\left(-s^2-4s-2\right)-2\kappa \left( \frac{2}{2+\mu t}-p\right)  \right]  \\
&=\frac{p^{N-1}}{4L(2+\mu t)} \left[(2s-s^3+4-2s^2)-(s^3+4s^2+2s)-2\kappa(3\mu t - \mu L^2t^3)\right]\\
&= \frac{p^{N-1}}{4L(2+\mu t)} \left[-2s^3-6s^2+4(1-Lt)+2Lt(L^2t^2-1)\right]\\
&= \frac{p^{N-1}}{2L(2+\mu t)} \left[-(s^3+3s^2+2s)+s(s+1)(s+2)\right] = 0,\\
S_{N,N}:&= \lambda_{N-1}A_{1,1}+\lambda_{N-1}'A'_{1,1}-\frac{1}{2\mu}\beta_{N}+\frac{1}{2L}\alpha_{N} \\
&= \frac{1}{4L} \cdot \frac{s+2}{2+\mu t} + \frac{1}{4L} \cdot \frac{s}{2+\mu t}-\frac{1}{2\mu}\left(1-\frac{s+2}{2+\mu t}+\frac{s}{2+\mu t}\right)\\
&=\frac{1}{2L}\cdot \frac{s+1}{2+\mu t} - \frac{1}{2\mu} \cdot \frac{\mu t}{2+\mu t}\\
&=\frac{1}{2}\cdot \frac{t}{2+\mu t} - \frac{1}{2} \cdot \frac{t}{2+\mu t} = 0.
\end{align*} 
\textbf{Diagonal Interior Entries.}  For the interior index $1\leq k \leq N-1$ we have,
\begin{align*}
S_{k,k}: &=A_{0,0}\lambda_{k}+A'_{0,0}\lambda'_{k} + A_{1,1}\lambda_{k-1} + A'_{1,1}\lambda'_{k-1} - \frac{1}{2\mu}\beta_{k}+\frac{1}{2L}\alpha_k\\
&= \frac{1}{4L}(\lambda_{k-1}+\lambda_{k-1}') + \frac{1}{4L}\left((2-s^2)\lambda_k+(-s^2-4s-2)\lambda'_{k}\right) - \frac{1}{2\mu}(\lambda_{k-1}'-\lambda_k' - (\lambda_{k-1}-\lambda_k)), \\
&= \left[\frac{(s+1)p}{2L(s+2)} + \frac{2-s^3-3s^2}{2L(s+2)} - \frac{1-p}{\mu(s+2)} \right]\lambda_k \\
&= \frac{1}{2L(s+2)}\left( (2+\mu t)\kappa p - 2\kappa + 2-s^3-3s^2 \right) \lambda_k \\
&= \frac{1}{2L(s+2)}\left( \kappa(\mu L^2 t^3 - 3\mu t + 2) - 2\kappa + 2-s^3-3s^2 \right) \lambda_k \\
& = \frac{1}{2L(s+2)}\left( L^3t^3 - 3Lt + 2 - s^3 -3s^2 \right) \lambda_k = 0.
\end{align*} 
This proves $\BS = 0$.
\end{proof}

\section{Lemmas for Theorem~\ref{thm:primal_dual_equal_case31}}
In this section we provide lemmas along with proofs that are needed for proving Theorem~\ref{thm:primal_dual_equal_case31}. In this regime, we assume $\Delta=1$, and we have
$$
\sqrt{3}/L< t < 2/L.
$$
We define
\begin{equation}\label{p3_proof:c_def}
c^{(N)} = \frac{s^{2N}\kappa}{1+s^{2N}(\kappa-1)} = \frac{(Lt-1)^{2N}L}{\mu+(Lt-1)^{2N}(L-\mu)}
\end{equation}
and
\begin{equation}
G^{(N)}_{i,j} = 2\mu c^{(N)}(-s)^{i+j-2N},~~\forall~0\leq j\leq i\leq N.
\end{equation}
\\
The following lemmas are used in the proof of Theorem~\ref{thm:primal_dual_equal_case31}. Lemma~\ref{r3:LemmaS} verifies the positive semi-definiteness of $\BS$ and Lemma~\ref{lem:r3c1:full-primal} shows that the primal solution $\BG$ satisfies the constraints \eqref{PP1} for all pairs, hence the optimal solution is also feasible for the full SDP~\eqref{def:sdp_pep}. 
\begin{lemma}\label{r3:LemmaS}
For fixed $N$ we have
$$
\BS = \sum_{k=0}^{N-1} \left(\lambda_k\BA_{k,k+1} + \lambda'_k \BA_{k+1,k}\right) + \diag \left(\frac{\balpha}{2L}-\frac{\bbeta}{2\mu}\right).
$$
where
\begin{align*}
\lambda_k&= \left(\frac{1}{s+1}+\frac{1}{s^{2k+1}(s+1)} \right)c,~\lambda'_k= \left(\frac{1}{s^{2k+1}(s+1)}-\frac{s}{s+1} \right)c,~~0\leq k\leq N-1, \\
\alpha_{k} &= 0,~~0\leq k\leq N,\\
\beta_k &= 0,~~0\leq k\leq N-1,~~\beta_{N}=1-c.    
\end{align*}
Let
\begin{equation}\label{def:T_s}
  T(s):=\frac{\log\!\left(\frac{2}{s+3}\right)}
     {2\log s}.
\end{equation}
We have 
\[
\BS\succeq0
\quad\Longleftrightarrow\quad
N\leq T(s)
\quad\Longleftrightarrow\quad
N\leq\lfloor T(s)\rfloor=N_0,
\].
\end{lemma}
\begin{proof}
\textbf{Weighted decomposition.} We prove the following decomposition of $\BS$:
\begin{equation}\label{r3p1:decompS}
\BS = \sum_{k=0}^{N-1} \omega_k (s \be_k + \be_{k+1})(s \be_k + \be_{k+1})^{\top} 
\end{equation}
where
$$
\omega_k = \frac{c((s+3)s^{2k+2}-2)}{4L(s+1)s^{2k+2}}.
$$
\textbf{Entrywise verification.} Now we verify by each diagonal and off-diagonal entry. For diagonal entries, we have
\begin{align*}
 S_{0,0}
 &=\frac{(2-s^2)\lambda_0-(s^2+4s+2)\lambda'_0}{4L}\\
 &=\frac{c}{4Ls}
       \left[(2-s^2)-(s^2+4s+2)(1-s)\right]\\
 &=\frac{c}{4L}(s^2+2s-2)
 =s^2 \omega_0.
\end{align*}
For an interior index $1\leq k\leq N-1$, there is no diagonal
contribution from $\balpha$ or $\bbeta$. Moreover,
\[
 \lambda_{k-1}+\lambda'_{k-1}
 =\frac{c}{s+1}\bigl(1-s+2s^{-2k+1}\bigr).
\]
Thus
\begin{align*}
 S_{k,k}
 &=\frac{\lambda_{k-1}+\lambda'_{k-1}
          +(2-s^2)\lambda_k-(s^2+4s+2)\lambda'_k}{4L}\\
 &=\frac{c}{4L(s+1)}
   \left[(1-s+2s^{-2k+1})
         +(2-s^2)(1+s^{-2k-1})\right.\\
 &\hspace{39mm}\left.
         -(s^2+4s+2)(s^{-2k-1}-s)\right]\\
 &=\frac{c}{4L(s+1)}
       \bigl((1+s^2)(s+3)-4s^{-2k}\bigr)\\
 &=\omega_{k-1}+s^2\omega_k.
\end{align*}
The terminal diagonal entry contains the multiplier
$\beta_N=1-c$. We obtain
\begin{align*}
 S_{N,N}
 &=\frac{\lambda_{N-1}+\lambda'_{N-1}}{4L}
       -\frac{1-c}{2\mu}\\
 &=\frac{c}{4L(s+1)}\bigl(1-s+2s^{-2N+1}\bigr)
       -\frac{c}{2L}(s^{-2N}-1)\\
 &=\frac{c}{4L(s+1)}
       \left[1-s+2s^{-2N+1}-2(s+1)(s^{-2N}-1)\right]\\
 &=\frac{c}{4L(s+1)}\bigl(s+3-2s^{-2N}\bigr)
 =\omega_{N-1}.
\end{align*}
For off-diagonal entries, we have
\begin{align*}
 S_{k,k+1}=S_{k+1,k}
 &=\frac{s\lambda_k-(s+2)\lambda'_k}{4L}\\
 &=\frac{c}{4L(s+1)}
   \left[s\bigl(1+s^{-2k-1}\bigr)
         -(s+2)\bigl(s^{-2k-1}-s\bigr)\right]\\
 &=\frac{c}{4L(s+1)}
       \bigl(s(s+3)-2s^{-2k-1}\bigr)
 =s \omega_k.
\end{align*}
\textbf{Positive-semidefiniteness threshold.} The preceding calculations establish
\eqref{r3p1:decompS}. In particular, for every
$\bx=(x_0,\ldots,x_N)^{\top}\in\mathbb R^{N+1}$,
\[
\bx^{\top}\BS\bx
=
\sum_{k=0}^{N-1}\omega_k(sx_k+x_{k+1})^2.
\]
Thus, $\BS\succeq0$ whenever $\omega_k\geq0$ for all
$0\leq k\leq N-1$.\\
Since $0<s<1$, $L>0$, and $c>0$, we have
\[
\operatorname{sign}(\omega_k)
=
\operatorname{sign}\bigl((s+3)s^{2k+2}-2\bigr).
\]
The sequence $(s+3)s^{2k+2}-2$ is strictly decreasing
in $k$. Thus
\[
\omega_k\geq0\quad\text{for all }0\leq k\leq N-1
\quad\Longleftrightarrow\quad
\omega_{N-1}\geq0.
\]
Moreover,
\[
\begin{aligned}
\omega_{N-1}\geq0
&\Longleftrightarrow
s^{2N}\geq\frac{2}{s+3}\\
&\Longleftrightarrow
2N\log s\geq\log\!\left(\frac{2}{s+3}\right)\\
&\Longleftrightarrow
N\leq
\frac{\log\!\left(\frac{2}{s+3}\right)}
     {2\log s}
= T(s),
\end{aligned}
\]
where the last inequality reverses direction because
$\log s<0$. This proves $\BS\succeq0$ for $N\leq T(s)$,
including equality.
Conversely, if $N>T(s)$, then $\omega_{N-1}<0$.
The decomposition gives
\[
\be_N^{\top}\BS\be_N
=S_{N,N}
=\omega_{N-1}<0.
\]
Therefore $\BS\not\succeq0$.\\
We conclude that
\[
\BS\succeq0
\quad\Longleftrightarrow\quad
N\leq T(s)
\quad\Longleftrightarrow\quad
N\leq\lfloor T(s)\rfloor=N_0,
\]
where the final equivalence uses that $N$ is an integer.
This completes the proof.
\end{proof}

\begin{lemma}
\label{lem:r3c1:full-primal}
We have
\[
 \lag\BA_{i,j},\BG\rag=\delta_i-\delta_j
 \qquad\text{for every }0\leq i\neq j\leq N.
\]
\end{lemma}

\begin{proof}
Fix $N\geq 1$ and suppress the superscript $(N)$ throughout the
proof. We calculate directly from the coefficient matrices in
\eqref{def:Aij} and~\eqref{def:AA}.\\
\smallskip
First we calculate the RHS. Since $G_{k,k+1}=-sG_{k,k}$ and
$\mu t(1-s)=\kappa^{-1}(1-s^2)$, the recursion for $\delta_k$
gives
\begin{equation}\label{eq:r3c1:full-gap-formula}
\begin{aligned}
 \delta_k
 &=c+\mu t(1-s)c\sum_{r=k}^{N-1}s^{2r-2N}\\
 &=c\left[1+\kappa^{-1}\bigl(s^{2k-2N}-1\bigr)\right],
 \qquad 0\leq k\leq N,
\end{aligned}
\end{equation}
where the sum is empty for $k=N$. In particular,
\begin{equation}\label{eq:r3c1:rankone-pair-identities}
 \delta_i-\delta_j
 =\frac{c}{\kappa}
   \bigl(s^{2i-2N}-s^{2j-2N}\bigr)
 =\frac{G_{i,i}-G_{j,j}}{2L}.
\end{equation} 
We reuse the computational result~\eqref{eq:rankone:pair-traces} since $\BG$ is again rank $1$. Fix $0\leq i<j\leq N$ with $h:=j-i\geq 1$, and set
$\rho:=-s$ and $H:=G_{i,i}>0$. We have
\begin{align*}
 \lag\BA_{i,j},\BG\rag
 &=\frac{G_{ii}}{4L} \left(1-\frac{L^2t^2}{(1+s)^2} \right) (1-(-s)^h)^2 + \frac{tG_{ii}}{2(1+s)}(1-s^{2h}) \\
 &= \frac{tG_{ii}}{2(1+s)}(1-s^{2h}) \\
 &=\frac{2c\mu t}{2Lt}(s^{2i-2N}-s^{2j-2N})\\
 &=\frac{c}{\kappa}
   \bigl(s^{2i-2N}-s^{2j-2N}\bigr) = \delta_i-\delta_j,
\end{align*}
where the first term vanishes since $(1+s)^2 = L^2t^2$ and hence
\begin{align*}
 \lag\BA_{j,i},\BG\rag
 &=\frac{G_{ii}}{4L} \left(1-\frac{L^2t^2}{(1+s)^2} \right) (1-(-s)^h)^2 - \frac{tG_{ii}}{2(1+s)}(1-s^{2h}) \\
 &= -\frac{tG_{ii}}{2(1+s)}(1-s^{2h}) =\delta_j-\delta_i.
\end{align*}
This verifies the computation for all the pairs.
\end{proof}

\end{document}